\documentclass[12pt, reqno]{amsart}
\usepackage{amssymb,mathrsfs,amsfonts,amsmath}
\usepackage[linktocpage=true]{hyperref}
\usepackage{enumerate}
\usepackage{dynkin-diagrams}

\begin{document}
\numberwithin{equation}{section}

\def\1#1{\overline{#1}}
\def\2#1{\widetilde{#1}}
\def\3#1{\widehat{#1}}
\def\4#1{\mathbb{#1}}
\def\5#1{\frak{#1}}
\def\6#1{{\mathcal{#1}}}

\def\C{{\4C}}
\def\R{{\4R}}
\def\N{{\4N}}
\def\Z{{\4Z}}

\title[Rationality of proper holomorphic maps between bounded symmetric domains]{Rationality of proper holomorphic maps between bounded symmetric domains of the first kind}
\author[S.-Y. Kim ]{Sung-Yeon Kim}
\address{S.-Y. Kim: Center for Complex Geometry, Institutue for Basic Science, 
55, Expo-ro, Yuseong-gu, Daejeon, Korea, 34126 }
\email{sykim8787@ibs.re.kr}


\maketitle


\def\Label#1{\label{#1}{\bf (#1)}~}


\def\cn{{\C^n}}
\def\cnn{{\C^{n'}}}
\def\ocn{\2{\C^n}}
\def\ocnn{\2{\C^{n'}}}


\def\dist{{\rm dist}}
\def\const{{\rm const}}
\def\rk{{\rm rank\,}}
\def\id{{\sf id}}
\def\aut{{\sf aut}}
\def\Aut{{\sf Aut}}
\def\CR{{\rm CR}}
\def\GL{{\sf GL}}
\def\Re{{\sf Re}\,}
\def\Im{{\sf Im}\,}
\def\span{\text{\rm span}}
\def\mult{\text{\rm mult\,}}
\def\reg{\text{\rm reg\,}}
\def\ord{\text{\rm ord\,}}
\def\hot{\text{\rm HOT\,}}

\def\codim{{\rm codim}}
\def\crd{\dim_{{\rm CR}}}
\def\crc{{\rm codim_{CR}}}

\def\eps{\varepsilon}
\def\d{\partial}
\def\a{\alpha}
\def\b{\beta}
\def\g{\gamma}
\def\G{\Gamma}
\def\D{\Delta}
\def\Om{\Omega}
\def\k{\kappa}
\def\l{\lambda}
\def\L{\Lambda}
\def\z{{\bar z}}
\def\w{{\bar w}}
\def\Z{{\1Z}}
\def\t{\tau}
\def\th{\theta}

\emergencystretch15pt
\frenchspacing

\newtheorem{Thm}{Theorem}[section]
\newtheorem{Cor}[Thm]{Corollary}
\newtheorem{Pro}[Thm]{Proposition}
\newtheorem{Lem}[Thm]{Lemma}

\theoremstyle{definition}\newtheorem{Def}[Thm]{Definition}

\theoremstyle{remark}
\newtheorem{Rem}[Thm]{Remark}
\newtheorem{Exa}[Thm]{Example}
\newtheorem{Exs}[Thm]{Examples}

\def\bl{\begin{Lem}}
\def\el{\end{Lem}}
\def\bp{\begin{Pro}}
\def\ep{\end{Pro}}
\def\bt{\begin{Thm}}
\def\et{\end{Thm}}
\def\bc{\begin{Cor}}
\def\ec{\end{Cor}}
\def\bd{\begin{Def}}
\def\ed{\end{Def}}
\def\br{\begin{Rem}}
\def\er{\end{Rem}}
\def\be{\begin{Exa}}
\def\ee{\end{Exa}}
\def\bpf{\begin{proof}}
\def\epf{\end{proof}}
\def\ben{\begin{enumerate}}
\def\een{\end{enumerate}}
\def\beq{\begin{equation}}
\def\eeq{\end{equation}}

\begin{abstract}
Let $D_{p,q}$ and $D_{p',q'}$ be irreducible bounded symmetric domains of the first kind with rank $q$ and $q'$, respectively and let $f:D_{p,q}\to D_{p',q'}$ be a proper holomorphic map that extends $C^2$ up to the boundary. 
In this paper we show that if $q, q'\geq 2$ and $f$ maps Shilov boundary of $D_{p,q}$ to Shilov boundary of $D_{p',q'}$, then $f$ is of the form
$f = \imath\circ F$, where
$$
	F=F_1\times F_2\colon D_{p,q}\to \Omega_1'\times \Omega_2', 
$$
$\Omega_1'$ and $\Omega_2'$ are bounded symmetric domains, $F_1 \colon D_{p,q}\to \Omega_1'$ is a proper rational map, $F_2:D_{p,q}\to \Omega_2'$ is not proper and $\imath: \Omega_1' \times \Omega'_2 \hookrightarrow D_{p',q'}$ is a holomorphic totally geodesic isometric embedding of a reducible bounded symmetric domain $\Omega_1' \times \Omega_2'$ into $D_{p',q'}$ with respect to canonical K\"ahler-Einstein metrics. Moreover, if $p>q$, then $f$ is a rational map.
As an application we show that a proper holomorphic map $f:D_{p,q}\to D_{p',q'}$ that extends $C^\infty$ up to the boundary is a rational map or a totally geodesic isometric embedding with respect to the Kobayashi metrics, if
$3\leq q \leq q'\leq 2q-1.$
\end{abstract}

%
%
%
%
%
%
%
%
%

\section{Introduction}
The purpose of this paper is to investigate the possibility of extending proper holomorphic maps between bounded symmetric domains to their compact duals, which are Hermitian symmetric manifolds of compact type.  
Hermitian symmetric manifolds of compact type are the cominuscle rational homogeneous projective varieties. In particular, they are birational to projective spaces. 
Therefore it is interesting to find conditions under which a proper holomorphic map between two bounded symmetric domains can be extended to a rational map between their compact duals.

Proper holomorphic maps from a unit disc to itself are the finite Blaschke products. In higher dimensional case, there are certain rigidity phenomena. These phenomena were first discovered by Poincar\'{e} (\cite{P07}), who proved that any biholomorphic map between two connected open
pieces of the unit sphere in $\mathbb C^2$ is a restriction of an automorphism of its compact dual. Later, Alexander (\cite{ A74}) and Henkin-Tumanov (\cite{TK82}) generalized
his result to higher dimensional unit balls and higher rank bounded symmetric domains
respectively.

For proper holomorphic maps between balls of different dimensions, Webster (\cite{Web79}) used the CR geometry of the sphere to show that a proper holomorphic map $f:\mathbb B^n
\to \mathbb B^{N}$ that is three times continuously differentiable up to the boundary is a restriction of projective linear embedding of $\mathbb P^n$ into $\mathbb P^{N}$ if $N=n+1$ and $n\geq 3.$ This result has been generalized by many mathematicians under certain conditions on the dimension difference and the regularity of $f$ on the boundary. We refer the works of Cima--Suffridge \cite{CS90}, Faran \cite{Fa86}, Ebenfelt \cite{Eb13}, Forstneri\v{c} \cite{F86, F89}, Globevnik \cite{G87}, Huang \cite{H99, H03}, Huang--Ji \cite{HJ01}, Huang--Ji--Xu \cite{HJX06}, Stens\o nes \cite{S96}, D'Angelo \cite{D88a, D88b, D91, D03}, D'Angelo--Kos--Riehl (\cite{DKR03}), D'Angelo--Lebl \cite{DL09, DL16} and the references therein.

When the dimension difference is arbitrary, Forstneri\v{c} \cite{F89} proved the rational extendability of $f:\mathbb B^n\to\mathbb B^N$.
More precisely, he proved:
\bt\label{Fo}
Let $U$ be an open ball centered at a point $P\in \partial\mathbb B^n$, and let $M = U\cap \mathbb B^n$.
If $N> n > 1$ and $f: \overline{\mathbb B}^n\cap U\to \mathbb C^N$ is a mapping of class $C^{N-n+ 1}$ that is holomorphic on
$\mathbb B^n\cap U$ and takes $M$ to the unit sphere $\partial\mathbb B^N$, then $f$ is rational, $f= (p_1,\ldots, p_N)/q$, where
the $p_j$ and $q$ are holomorphic polynomials of degree at most $N^2(N-n+ I)$. The
extended map is holomorphic on $\mathbb B^n$, it maps $\mathbb B^n$ to $\mathbb B^N$, and it has no poles on $\partial \mathbb B^n$.
\et

In contrast to the case of complex unit balls, which are
precisely the bounded symmetric domains of rank 1, much less is known for proper holomorphic maps between bounded symmetric domains $\Omega, \Omega'$ of higher rank. 
By using the geometric properties of Hermitian symmetric spaces with respect to the canonical K\"{a}hler-Einstein metrics, Tsai (\cite{Ts93}) proved that a proper holomorphic map $f:\Omega\to\Omega'$
is a restriction of a totally geodesic isometric embedding of the compact dual, if $\text{rank}(\Omega)\geq \textup{rank}(\Omega')\geq 2$. 
When $\text{rank}(\Omega)< \textup{rank}(\Omega')$, total geodesy of $f$ fails in general. 
But when the rank difference or dimension difference is sufficiently small, classification or nonexistence of proper holomorphic
maps were obtained for certain pairs of irreducible bounded symmetric domains of rank $\geq 2$. We refer the readers to Chan \cite{C20, C21}, Henkin-Novikov \cite{HN84}, Kim-Mok-Seo \cite{KMS22}, Kim-Zaitsev \cite{KZ13, KZ15}, Mok \cite{M08c}, Mok-Ng-Tu \cite{MNT10}, Ng \cite{N13, N15a, N15b}, Seo \cite{S15, S16, S18} and Tu \cite{Tu02a, Tu02b}. 

As for the rational extendability of $f$, it is shown in \cite{KMS22} that if $\Omega$ and $ \Omega'$ are of the same type or of type three and type one, respectively and if the rank difference is sufficiently small($2\leq \textup{rank}(\Omega')<2~\textup{rank}(\Omega)-1$),
then $f$ is of the form
$f = \imath\circ F$, where
$$
	F=F_1\times F_2\colon \Omega\to \Omega_1'\times \Omega_2', 
$$
$\Omega_1'$ and $\Omega_2'$ are bounded symmetric domains, $F_1 \colon \Omega\to \Omega_1'$ is a totally geodesic isometric embedding and $\imath: \Omega_1' \times \Omega'_2 \hookrightarrow \Omega'$ is a holomorphic totally geodesic isometric embedding of a reducible bounded symmetric domain $\Omega_1' \times \Omega_2'$ into $\Omega'$ with respect to canonical K\"ahler-Einstein metrics. As a consequence $f$ has a factor $F_1$ that extends to a totally geodesic isometric embedding of the compact dual of $\Omega$ into the compact dual of $\Omega_1'$.
If $\Omega$ is of rank one, there are proper holomorphic maps which do not have the above property. Reiter-Son \cite{RS} and Xiao-Yuan \cite{XY20} classified proper holomorphic maps from $\mathbb B^n$ to type four bounded symmetric domains. Those maps are algebraic but not necessarily have the above property. 
\medskip

In this paper, we generalize the result of Forstneri\v{c}. More precisely, we show the rational extendability of proper holomorphic maps $f$ between bounded symmetric domains of type one with rank $\geq 2$ under a certain condition on the boundary values. If a proper holomorphic map $f$ extends continuously to the boundary, $f$ maps the boundary of $\Omega$ into the boundary of $\Omega'$. The topological boundary of $\Omega$ is a disjoint union of the $G$-orbits $S_r(\Omega)$, $r=0,\ldots,\text{rank}(\Omega)-1$, where $G$ is the identity component of $\Aut(\Omega)$ such that each $S_r(\Omega)$ is foliated by bounded symmetric domains of rank $r$. Among those orbits, Shilov boundary $S_0(\Omega)$ is the unique closed orbit and contains no complex manifold of positive dimension. In this paper, we impose a condition that $f$ maps Shilov boundary of $\Omega$ to Shilov boundary of $\Omega'$, which is a generalization of the condition that $f:\overline{\mathbb B}^n\cap U\to \mathbb C^N$ maps $\partial\mathbb B^n\cap U$ to $\partial\mathbb B^N$.

\begin{Thm}\label{main}
Let $D_{p,q}$, $D_{p',q'}$ be irreducible bounded symmetric domains of type one with rank $q$, $q'$, respectively and let $f:D_{p,q}\to D_{p',q'}$ be a proper holomorphic map.
Suppose there exist a point $P\in S_0(D_{p,q})$ and an open neighborhood $U$ of $P$ such that $f$ extends $C^2$ to $ U$. Suppose further that $q, q'\geq 2$ and $f$ maps $S_0(D_{p,q})\cap U$ to $S_0(D_{p',q'})$.
Then $f$ is of the form
$f = \imath\circ F$, where
$$
	F=F_1\times F_2\colon D_{p,q}\to \Omega_1'\times \Omega_2', 
$$
$\Omega_1'$ and $\Omega_2'$ are bounded symmetric domains, $F_1 \colon D_{p,q}\to \Omega_1'$ is a proper rational map, $F_2:D_{p,q}\to \Omega_2'$ is not proper and $\imath: \Omega_1' \times \Omega'_2 \hookrightarrow D_{p',q'}$ is a holomorphic totally geodesic isometric embedding of a reducible bounded symmetric domain $\Omega_1' \times \Omega_2'$ into $D_{p',q'}$ with respect to canonical K\"ahler-Einstein metrics. Moreover, if $p>q$, then $f$ is a rational map.
\end{Thm}

We remark that if $q=1$, i.e., $D_{p,q}=\mathbb B^p$ and if $f$ maps Shilov boundary to Shilov boundary, then $f$ is a proper holomorphic map from $\mathbb B^p$ to $\mathbb B^{\dim D_{p',q'}}$. Therefore by Theorem \ref{Fo}, $f$ is a rational map if $f$ is sufficiently smooth up to the boundary and $p>1$. If $q'=1$, i.e., $D_{p',q'}=\mathbb B^{p'}$, then there is no proper holomorphic map $f:D_{p,q}\to \mathbb B^{p'}$ that extends continuously to an open neighborhood of a boundary point if $q>1$. 
\medskip

If $F_1$ in Theorem \ref{main} is a totally geodesic isometric embedding with respect to the canonical K\"{a}hler-Einstein metrics, then by \cite{M22}, $f$ is a totally geodesic isometric embedding with respect to Kobayashi metric. With this result and \cite{K21}, we obtain the following corollary.


\bc\label{cor}
Let $f:D_{p,q}\to D_{p',q'}$ be a proper holomorphic map that extends $C^\infty$ up to the boundary. Suppose 
$$3\leq q \leq q'\leq 2q-1.$$
Then $f$ is a rational map or of the form
$f = \imath\circ F$, where
$$
	F=F_1\times F_2\colon D_{p,q}\to \Omega_1'\times \Omega_2', 
$$
$\Omega_1'$ and $\Omega_2'$ are bounded symmetric domains, $F_1 \colon D_{p,q}\to \Omega_1'$ is a proper rational map, $F_2:D_{p,q}\to \Omega_2'$ is not proper and $\imath: \Omega_1' \times \Omega'_2 \hookrightarrow D_{p',q'}$ is a holomorphic totally geodesic isometric embedding of a reducible bounded symmetric domain $\Omega_1' \times \Omega_2'$ into $D_{p',q'}$ with respect to canonical K\"ahler-Einstein metrics. As a consequence $f$ extends rationally to the compact dual or a totally geodesic isometric embedding with respect to the Kobayashi metrics. 
\ec

Our method is to combine K\"{a}hler geometry of Hermitian symmetric spaces and CR geometry on the boundary orbits of bounded symmetric domains. 
We first generalize a strategy of characteristic bundles over a Hermitian symmetric space and moduli maps induced by $f$,
which was first used in the work of Mok-Tsai (\cite{MT92}). Using the properness of $f$, we define a moduli map $f^\sharp_r$ between the moduli spaces of invariantly totally geodesic subdomains of $\Omega$ and $\Omega'$, which is meromorphic. As in \cite{MT92}, the pseudoconcavity of the moduli space of invariantly totally geodesic subdomains forces the moduli map to extend globally to a rational map between moduli spaces of invariantly totally geodesic subspaces of the compact duals (cf. \cite{KMS22}). 

For each point $P$ in the compact dual $X$ of $\Omega$, we define a complex variety $\mathscr Z_P^r\subset X$ by the union of all characteristic subspaces of rank $r$ passing through $P$. If $P$ is a point in the Shilov boundary, then we can define a subset $\mathscr S_P^r\subset \mathscr Z_P^r$ by the union of all boundary components of rank $r$ with $P$ in their closure.  
We show that $\mathscr S_P^r$ is a CR submanifold in $S_r(\Omega)$ and $\mathscr Z_P^r$ is the smallest compact complex variety that contains $\mathscr S_P^r$. Moreover, if $f$ maps Shilov boundary to Shilov boundary, then the image of $\mathscr S_P^r$ under $f$ is determined by the second jet of $f$ restricted to the rank $r$ boundary orbit $S_r(\Omega)$. We use the CR geometry on the boundary orbits of $\Omega$ and $\Omega'$ to analyse the CR second fundamental form of $f$ and then show that $f$ decomposes into a product of two holomorphic maps, where one of the factors is extended rationally to a neighborhood of $P$ via lifting and pushing down the moduli map $f^\sharp_r$ through double fibration of the universal family of invariantly totally geodesic subspaces.

The organization of the current article is as follows. In Section 2, we describe the moduli spaces and universal family of invariantly totally geodesic subspaces for Hermitian symmetric space of type one.
In Section 3, we define moduli maps induced by $f$. Then we obtain a condition for the moduli maps to be lifted to the universal family of invariantly totally geodesic subspaces.
In Section 4, we investigate the CR structures of $S_r(\Omega)$ and $\mathscr S_P^r$. Then, in Section 5, we define CR second fundamental form of $f$, which 
describes the image of $\mathscr S_P^r$ via $f$. Throughout Section 4 and 5, we use the Einstein summation convention unless stated otherwise. Finally, in Section 6, we prove Theorem \ref{main} and Corollary \ref{cor}. 
\medskip

{\bf Acknowledgement} Author was supported by the Institute for Basic Science (IBS-R032-D1-2021-a00).

\section{Preliminaries}
In this section, we present some basic notions of bounded symmetric domains of type one. 
Then we will prove some basic properties of them. We refer \cite{KMS22}, \cite{Mok86} and \cite{MT92} as references.
\subsection{Hermitian symmetric spaces of type one}
For positive integers $p\geq q\geq 1$, define a basic form $\langle~,~ \rangle=\langle~,~\rangle_{p,q}$ on $\mathbb C^{p+q}$
by
$$\langle z, w\rangle=-\sum_{1\leq j\leq q}z_j\overline w_j+\sum_{q+1\leq j \leq p+q}z_j\overline w_j. $$
The noncompact Hermitian symmetric space $D_{p,q}$ of type one is the set of all $q$-planes $V \subset \mathbb C^{p+q}$ such that the
restriction $\langle~, ~\rangle\big|_V$ is negative definite. 
In Harish-Chandra coordinates, $D_{p,q}$ in the complex Grassmannian $Gr(p,q)$ is realized as follows:

Write $M^{\mathbb C}(p,q)$ for the set of $p \times q$ matrices with coefficients in $\mathbb C$, and denote by $\{e_1, \ldots, e_{p+q}\}$ the standard basis of $\mathbb C^{p+q}$.  For $Z \in M^\mathbb C(p,q)$, denoting by $v_k$, $1 \le k \le q$, the $k$-th column vector of $Z$ as a vector in $\mathbb C^p = {\rm Span}_\mathbb C\{e_{1+q},\ldots,e_{p+q}\}$ we identify $Z$ with the $q$-plane in $\mathbb C^{p+q}$ spanned by $\{e_k + v_k: 1 \le k \le q\}$.
Then we have
\begin{equation*}\label{typeI}
D_{p,q}= \left\{ Z\in M^{\mathbb C}(p,q): I_q -  Z^*Z >0 \right\},
\end{equation*}
where $Z^*$ denotes the conjugate transpose of $Z$.
Throughout the paper, we always assume that $D_{p,q}$ is defined as above.

The topological boundary $\partial D_{p,q}$ of $D_{p,q}$ is a disjoint union of the boundary orbits
$S_{p,q,r}$, $r=0,\ldots,q-1$, where each $S_{p,q,r}$ consists of all $q$-planes $V\subset \mathbb C^{p+q}$ such that the restriction  $\langle~, ~\rangle\big|_V$
has $r$ negative
and $q-r$ zero eigenvalues. The connected identity component $G$ of the biholomorphic automorphism
group $\Aut (D_{p,q} )$ is identified with the identity component of the group of all linear transformations of $\mathbb C^{p+q}$ preserving $\langle~, ~\rangle$. Therefore each $S_{p,q,r}$ is a $G$-orbit in $Gr(p,q)$. The rank $q-1$ boundary orbit $S_{p,q,q-1}$ is called the hypersurface boundary and the rank $0$ boundary orbit $S_{p,q,0}$ is called \emph{Shilov boundary}. We remark that $S_{p,q,0}$ is the set of extremal points of $\overline {D}_{p,q}$.

For a pair of linear subspaces $(V_1, V_2)$ such that $V_1\subset V_2$, denote by $[V_1, V_2]_q$ the set of all elements $x\in Gr(p,q)$ such that 
\begin{equation}\label{subgra1}
V_1 \subset V_x \subset V_2,
\end{equation}
where $V_x$ is the $q$-dimensional subspace of $\mathbb C^{p+q}$ corresponding to $x\in Gr(p,q)$.  
For a linear subspace $A\subset \mathbb C^{p+q}$, denote by $A^*$ the space of all vectors $v$ in $\mathbb C^{p+q}$ such that
$$\langle \cdot, v\rangle\big|_{A}=0.$$
$A$ is called a \emph{null space} if 
$$\langle A, A\rangle=0$$
or equivalently,
$$A\subset A^*.$$
If $A$ is a $(q-r)$-dimensional null space,
then $A^*$ is a $(p+r)$-dimensional subspace containing $A$ and
$[A, A^*]_q\cap S_{p,q,r}$ is a maximal complex submanifold in $\partial D_{p,q}$ biholomorphic to $D_{p-q+r, r}$, called a \emph{boundary component} of  rank $r$. Up to the action of $\Aut(D_{p,q})$, a boundary component of rank $r$ is
equivalent to a standard boundary component 
$$
F_r=\left\{\begin{pmatrix}
I_{q-r}&0\\
0&Z\\
\end{pmatrix}
, Z\in D_{p-q+r, r}\right\}.
$$
We remark that $[A, A^*]_q$ is the compact dual of $[A, A^*]_q\cap S_{p,q,r}$. Conversely, every maximal complex submanifold in $\partial D_{p,q}$ is of this form. See \cite{W72}.

%

%
%
%

\subsection{Invariantly totally geodesic subspaces}\label{moduli spaces}
An invariantly totally geodesic subspace of subdiagram type in $X=Gr(p,q)$ is a subgrassmannian $[V_1, V_2]_q$ for some linear subspaces $V_1$ and $V_2$ such that $V_1\subset V_2$.
Hence, for fixed positive integers $a \leq  b$, the moduli space of invariantly totally geodesic subspaces of subdiagram type with $\dim V_1 =a$, $\dim V_2 = b$ is the flag variety 
\begin{equation}
\mathcal F_{a,b}(X)=\{ (V_1, V_2) : \{0\} \subset V_1 \subset V_2\subset \mathbb C^{p+q}, \dim V_1 = a, \, \dim V_2 = b \}.
\end{equation}
For $\sigma=(V_1, V_2)\in \mathcal F_{a,b}(X)$ we sometimes denote the corresponding
subgrassmannian $[V_1, V_2]_q$ by $X_\sigma$.

A subgrassmannian $[V_1, V_2]_q\subset Gr(p,q)$ with $\dim V_1=q-r$ and $\dim V_2=p+r$ is called a \emph{characteristic subspace} of rank $r$. We denote the moduli space of characteristic subspaces of rank $r$ for $r=1,\ldots, q-1$ by $\mathcal D_r(X)$, i.e.,
\begin{equation}\label{css1}
\mathcal D_r(X)=\{ (V_1, V_2) : \{0\} \subset V_1 \subset V_2\subset \mathbb C^{p+q}, \dim V_1 = q-r, \, \dim V_2 = p+r \}.
\end{equation}
Then $\mathcal D_r(X)$ is biholomorphic to $G/P_r$ for a parabolic subgroup $P_r$ of $G$ and 
the automorphism group is $SL(p+q,\mathbb C)$ for $r> 0$ (see section 3.3 in \cite{Akhiezer}).
In particular, $\mathcal D_r(X)$ is a rational homogeneous manifold.
\medskip

For a given $\Omega=D_{p,q}$ and its compact dual $X=Gr(p,q)$, define
$$
	\mathcal{D}_r(\Omega):=\{\sigma\in \mathcal{D}_r(X)\colon  \Omega_\sigma := X_\sigma\cap\Omega\neq \emptyset\}.$$
Here $\Omega_\sigma$ is a bounded symmetric domain of type one with rank $r$. We call $\Omega_\sigma$ a characteristic subdomain of rank $r$.
For a rank $k$ boundary orbit $S_k=S_k(\Omega)=S_{p,q,k}$ of $\Omega$ with $k\geq r$, define
$$
	\mathcal{D}_r(S_k):=\{\sigma\in \mathcal{D}_r(X)
	\colon \Omega_\sigma:= X_\sigma\cap S_k \text{ is open in }X_\sigma\}.
$$
Then $\mathcal D_r(\Omega)$ and $\mathcal{D}_r(S_k)$, $k=r,\ldots,q-1$ are $G$-orbits in $\mathcal{D}_r(X)$ such that $\mathcal{D}_r(S_k)\subset \partial \mathcal{D}_r(\Omega)$. Furthermore $\mathcal D_r(S_r)$ is the unique closed orbit of $G$ such that $\sigma\in \mathcal D_r(S_r)$ if and only if $\Omega_\sigma$ is a boundary component of rank $r$. 

\subsection{Associated characteristic bundle}\label{cb} 
A unit vector $v\in T_xX$ is called a \emph{characteristic vector} if it realizes the maximum of the holomorphic sectional curvature of $X$ with respect to the canonical K\"{a}hler-Einstein metric. Characteristic fiber bundle $\mathscr C(X)$ is the bundle
$$\mathscr C(X)=\bigcup_{x\in X}\{[v]\in \mathbb PT_xX: v \text{ is a characteristic vector}\}$$
which is a holomorphic fiber bundle over $X$. When a Hermitian symmetric domain is realized by a bounded symmetric domain $\Omega$ via Harish-Chandra embedding in $X$, the fibers $\mathscr C_x$ over $x\in \Omega$ are
parallel with respect to Harish-Chandra coordinates.

For each characteristic vector $v\in T_xX$, there is an orthogonal
decomposition
$$ T_xX = \mathbb C v\oplus \mathscr H_v\oplus \mathscr N_v$$
into eigenspaces of the Hermitian form
$$\mathcal R_v(\xi, \eta):=R_{v\bar v\xi\bar\eta}$$
corresponding to eigenvalues $R_{v\bar v v\bar v}, 1/2 R_{v\bar v v\bar v}$ and $0$, where $R$ is the curvature tensor of the canonical K\"{a}hler-Einstein metric on $X$. 
Let $\mathbb P_v$ be the projective line in $X$ tangent to $v$ that realizes the maximum of holomorphic sectional curvature. Then there exists a unique characteristic subspace $X_\sigma$ of rank $q-1$ passing through $x$ such that the tangent space of $X_\sigma$ at $x$ is given by
$$T_x X_\sigma=\mathscr N_v.$$
Conversely, for each characteristic subspace $X_\sigma$ of rank $q-1$ passing through $x$, up to a constant multiple, there exists a unique characteristic vector $v_x$ at $x$ such that 
$$T_xX_\sigma=\mathscr N_{v_x}.$$
Moreover, for $x, y\in X_\sigma\cap \overline\Omega$, $v_x$ and $v_y$ are parallel in Harish-Chandra coordinates.
In either case, $\mathbb P_v\times X_\sigma$ is a totally geodesic submanifold in $X$ for $v=v_x$ or $v=v_y$.
Define
$$\mathscr {N}_{q-1}(X)=\bigcup_{[v]\in \mathscr C(X)}\mathscr N_v=\{T_xX_\sigma: \sigma\in \mathcal D_{q-1}(X), x\in X_\sigma\}.$$
It is a complex homogeneous fiber bundle over $X$, called the \emph{associated characteristic bundle} of rank $q-1$. 
Likewise we can define associated characteristic bundles of rank $r$ by
$$\mathscr {N}_r(X)=\{T_xX_\sigma: \sigma\in \mathcal D_r(X), x\in X_\sigma\}, \quad r=1,\ldots q-1,$$
which are complex homogeneous fiber bundles over $ X$.

Let $A\subset \mathscr C_x(X)$ be a set of characteristic directions at $x$.
Denote
$$\mathscr N_A=\bigcap_{[v]\in A}\mathscr N_v.$$
Then there exists a unique invariantly totally geodesic subspace $X_A$ passing through $x$ such that 
$$T_xX_A=\mathscr N_A.$$
Furthermore, there exists a unique invariantly totally geodesic subspace $Y_A$ of subdiagram type passing through $x$ such that 
$$\mathscr N_A=\mathscr N_{\mathcal C_x(X)\bigcap\mathbb P T_xY_A}.$$
Let 
$$\mathscr C_x(X_A):=\{[v]\in \mathscr C_x(X): v\in T_xX_A\}.$$
Then
$$T_xY_A=\mathscr N_{\mathscr C_x(X_A)}.$$
We will denote by $\mathscr H_A$ the orthogonal complement of $T_xY_A+ \mathscr N_A.$
Then we have an orthogonal decomposition
\beq\label{o-decomp}
T_xX=T_xY_A\oplus \mathscr H_A\oplus \mathscr N_A.
\eeq

\subsection{Universal family of characteristic subspaces}\label{am}
For a given $r>0$, define
$$\mathcal{U}_{r}(X):=\{(x, \sigma)\in X\times \mathcal{D}_{r}(X)\colon x\in X_\sigma\}.$$
Then there is a canonical double fibration 
$$
	\rho_{r}\colon \mathcal{U}_r(X)\to \mathcal{D}_r(X),
	\quad
	\pi_r\colon\mathcal{U}_r(X)\to X
$$
given by
$$\rho_r(x, \sigma)=\sigma,\quad \pi_r(x,\sigma)=x.$$
Note that
$$X_\sigma=\pi_r\left(\rho^{-1}_r(\sigma)\right).$$
Define 
$\imath_r\colon \mathcal{U}_r(X)\to \mathcal{G}(n_r, TX)$ with $n_r=\dim X_\sigma$ by $\imath_r(x,\sigma)=T_x X_\sigma$, where $ \mathcal{G}(n_r, TX)$ is a Grassmannian bundle over $TX$.
Then $\imath_r$ is a $G$-equivariant holomorphic embedding. Hence we may regard $\mathcal U_r(X)$ as a complex manifold in $\mathcal G(n_r, TX)$. For each $x\in X$, define
$$\mathcal Z_x^r:=\rho_r(\pi_r^{-1}(x)).$$ 
Similarly, we can define a $G$-equivariant holomorphic embedding $j_r\colon \mathcal{U}_r(X)\to \mathcal{G}(m_r, T\mathcal D_r(X))$ with $m_r=\dim \mathcal Z_x^r$ by $j_r(x,\sigma)=T_\sigma \mathcal Z_x^r$ and we may regard $\mathcal U_r(X)$ as a complex manifold in $\mathcal G(m_r, T\mathcal D_r(X))$(See \cite{HoN21}).

For a complex submanifold $M\subset X$, define
$$
	\mathcal Z_M^r:=\{\sigma\in \mathcal D_r(X): M\subset X_\sigma\},
\quad
\mathcal S_M^r:=\mathcal Z_M^r\cap \mathcal D_r(S_r),
$$
$$
	\mathscr Z_M^r:=\pi_r\left (\rho_r^{-1}(\mathcal Z_M^r)\right),
	\quad 
	\mathscr S_M^{r}:=\mathscr Z_M^r\cap S_r.
$$
For $M=X_\sigma$ and $M=\{x\}$, we will denote $\mathcal Z_M^r$ by $\mathcal Z_\sigma^r$ and $\mathcal Z_x^r$, respectively for simplicity.
We remark that If $M=[A, B]_q$, then
$$ \mathcal Z_M^r=\{\sigma=(V, W)\in \mathcal D_r(X): V\subset A,~B\subset W\}$$
and
\beq\label{ZM}
\mathscr Z_M^r=\{x\in X: \dim V_x\cap A\geq q-r,~\dim (V_x+ B)\leq p+r    \}.
\eeq
For a given $r$, we will omit superscript $r$ if there is no confusion.

\bl\label{subvar}
Let $M\subset X$ be a subgrassmannian. Then $\mathcal Z_M^r$ is a projective algebraic manifold covered by a finite union of Euclidean coordinate charts and
$\mathscr Z_M^r$ is a complex variety in $X$. Moreover, if $M$ is a compact dual of a rank $s$ boundary component, i.e. $M=X_\sigma$ with $\sigma\in \mathcal D_s(S_s)$, then $\mathcal S_M^r$ and $\mathscr S_M^r$ are smooth manifolds.
\el

\bpf
Let 
$M=[V_M, W_M]_q$ be a subgrassmannian in $X$. Then
$\sigma=(A, B)\in \mathcal D_r(X)$ is contained in $\mathcal Z_M^r$ if and only if
$$A\subset V_M,\quad W_M\subset B.$$
Therefore $\mathcal Z_M^r$ is biholomorphic to $Gr(q-r, V_M)\times Gr(r-s, \mathbb C^{p+q}/W_M)$.
Since $\pi_r$ is a proper holomorphic map, $\mathscr Z_M^r=\pi_r(\rho^{-1}_r(\mathcal Z_M^r))$ is a complex analytic variety in $X$.

Now suppose $M$ is of the form
$$M=[V_M, V_M^*]_q$$
for some $(q-s)$-dimensional null space $V_M$. Then
$$\mathcal S_M^r=\{(A, A^*):A\in Gr(q-r, V_M)\}$$
and
$$\mathscr S_M^r=\{x\in S_r: \dim V_x\cap V_M=q-r\},$$
which completes the proof. 
\epf

\bl
Suppose $\Omega$ is of tube type, i.e., $\Omega=D_{q,q}.$
Then for a Shilov boundary point $P\in S_0$, 
$\mathcal S_{P}^r$ is a maximal totally real submanifold of $\mathcal Z_P^r$.
\el

\bpf
Let $\sigma=(A, B)\in \mathcal Z_P^r$. Then $A$ is a subspace of $V_P$, implying that
$$\langle A, A\rangle=0$$
and
$\sigma\in \mathcal D_r(S_r)$ if and only if 
$$ B=A^*.$$
Hence under a biholomorphic equivalence $\mathcal Z_P\simeq Gr(q-r, V_P)\times Gr(q-r, \mathbb C^{2q}/V_P)$, $\mathcal S_P$ is identified with
$$\{(A, A^*/V_P):A\in Gr(q-r, V_P)\}\subset Gr(q-r, V_P)\times  Gr(q-r, \mathbb C^{2q}/V_P) .$$
Therefore there exists a totally real embedding
$\imath:Gr(q-r,  V_P)\to \mathcal Z_P$ defined by
$\imath(A)=(A, A^*)$ whose image coincides with $ \mathcal S_P$, which completes the proof.
\epf


We regard $X=Gr(p,q)$ as a submanifold in a projective space $\mathbb P^N$ via the first canonical embedding $\imath:X\to \mathbb P^N$.
Then for each invariantly totally geodesic subspace $X_\sigma$ of subdiagram type, there exists a unique $\ell$-dimensional projective space $\mathbb P^\ell_\sigma\subset \mathbb P^N$ such that
$$\imath(X_\sigma)=\imath(X)\cap \mathbb P^\ell_\sigma.$$
$\mathbb P^\ell_\sigma$ is the smallest projective space in $\mathbb P^N$ that contains $\imath(X_\sigma)$.  
Therefore the manifold
$$ \mathcal{U}_{a,b}(X):=\{(x, \sigma)\in X\times \mathcal F(a,b;\mathbb C^{p+q}):x\in X_\sigma\}$$
can be regarded as a submanifold in $\mathbb P^N\times \mathscr G(\ell, N)$, where $\mathscr G(\ell, N)$ is the set of all $\ell$-dimensional projective spaces  in $\mathbb P^N$. 

Let 
$$\mathcal M=\{(x, L)\in \mathbb P^N\times \mathscr G(\ell, N): x\in L\}$$
be the tautological $\mathbb P^\ell$-bundle over $\mathscr G(\ell, N)$.
Then $\mathcal M$ is defined locally by a quadratic equation
$$Q(x,L)=(Q_1,\ldots, Q_{N-\ell})(x,L)=0$$
with the property
$$\{ x\in \mathbb P^N:Q(x,L)=0\}=L.$$
For a point $P\in \mathbb P^N$, define
$$\mathcal M_P:=\{P\}\times \{L\in \mathscr G(\ell, N): P\in L\}.$$
Then $\mathcal M_P$ is a projective algebraic manifold defined locally by
$$Q(P,\cdot)=0.$$
Choose an open set $O\subset \mathbb P^N\times \mathscr G(\ell, N)$ on which
$Q=0$ can be expressed by
$$Q_j(x, L)=\sum_{k=0}^N a_{j,k}x_k=0,\quad j=1,\ldots, N-\ell,$$
where $x=[x_0;x_1;\cdots;x_N]\in \mathbb P^N$ and $(a_{j,k})_{j,k}$ is a matrix of size $(N-\ell)\times(N+1)$.
We may assume $x_0=1$ and $(a_{j,k})_{j,k=1}^{N-\ell}$ is an identity matrix.
A complex manifold $\mathcal N\subset \mathcal M_P\cap O$ is defined locally by 
$\{(P, h(\zeta)) \in \mathcal M_P: \zeta\in U\}$ for some holomorphic embedding $h:U\to \mathscr G(\ell, N)$, where $U$ is a connected open set in a complex Euclidean space.
Then $h(\zeta)=(h_{j,k}(\zeta))$ 
satisfies
$$A_j(P;\zeta):=Q_j(P, h(\zeta))=\sum_{k=0}^N h_{j,k}(\zeta)P_k=0,\quad j=1,\ldots,N-\ell.$$
By taking the derivatives with respect to $\zeta$ at $\zeta_0\in U$, we obtain a system of affine equations 
$A^{(m)}(\cdot;\zeta_0)=((\partial_\zeta^\alpha A_j)(\cdot;\zeta_0): j=1,\ldots,N-\ell, ~|\alpha|\leq m)$
defined by
$$(\partial_\zeta^\alpha A_j)(P;\zeta_0)=\sum_{k=0}^N (\partial_\zeta^\alpha h_{j,k})(\zeta_0) P_k,\quad j=1,\ldots, N-\ell$$
whose coefficients depend holomorphically on $\zeta_0$, where for $a\in \mathbb N^d$ with $d=\dim U$,
$$\partial^\a_\zeta=\partial_{\zeta_1}^{\a_1}\cdots\partial_{\zeta_d}^{\a_d}$$
and 
$$|\a|=\a_1+\cdots+\a_d.$$ 
By linear algebra, we obtain the following lemma.

\bl\label{max rank}
Let $\mathcal N\subset \mathcal M_P\cap O$ be a complex manifold. Then
$$\bigcap_{L\in \mathcal N}L=\{P\}$$
if and only if there exists an integer $m$ and a point $(P,h(\zeta_0))\in \mathcal N$ such that 
$A^{(m)}(\cdot;\zeta_0)$ is of maximal rank. 
\el


\section{Induced moduli maps}
We identify $X'=Gr(p',q')$ with its image in a projective space $\mathbb P^{N'}$ via the first canonical embedding.
For each $r>0$, there exists a pair of integers $(a,b)=(a_r, b_r)$ depending only on $r$ such that for general $\sigma\in \mathcal D_r(\Omega)$, the smallest subgrassmannian of $X'$ that contains $ f(\Omega_\sigma)$ is of the form $[A_\sigma, B_\sigma]_{q'}$ with $\dim A_\sigma=a,~\dim B_\sigma=b$.  
Define a map
$f^\sharp_r:\mathcal D_r(\Omega)\to \mathcal F_{a,b}(X)$
by 
$$f^\sharp_r (\sigma)=(A_\sigma, B_\sigma).$$
Then as in \cite{KMS22}, $f^\sharp_r$ extends to a rational map $f^\sharp_r:\mathcal D_r(X)\to \mathcal F_{a,b}(X')$.
Define  
$$ \mathcal U'_{a,b}:=\{(y, L)\in X'\times \mathcal F_{a,b}(X'): y\in L\}.$$
Then 
there exists a canonical double fibration 
$$\rho'_{a,b}\colon \mathcal{U}_{a,b}'\to \mathcal F_{a,b}(X'),
\quad
\pi_{a,b}'\colon\mathcal{U}'_{a,b}\to X'. $$
The map $\mathcal F_r\colon \mathcal U_r(\Omega)\to \mathcal U'_{a,b}$ defined by $\mathcal F_r(x, \sigma)=(f(x), f^\sharp_r(\sigma))$ preserves the double fibrations, i.e.
$$\pi'_{a,b}\circ \mathcal F_r=f\circ \pi_r,\quad \rho'_{a,b}\circ\mathcal F_r=f^\sharp_r\circ\rho_r,$$
where 
$$\rho_r:\mathcal U_r(X)\to \mathcal D_r(X),\quad \pi_r:\mathcal U_r(X)\to X$$
is the canonical double fibration of the universal characteristic bundle over $X$.

We will regard $\mathcal U'_{a,b}$ as a submanifold of $\mathbb P^{N'}\times \mathscr G(\ell, N')$ for some $\ell$ and $N'$.
Similar to the case of universal characteristic bundles, for a complex submanifold $M\subset X'$ we define
$${{\mathcal Z}_M'}^{a,b}:=\{(A, B)\in \mathcal F_{a,b}(X'): M\subset [A, B]_{q'}\}$$
and 
$${\mathscr {Z}_M'}^{a,b}=\pi_{a,b}'\left(\left(\rho_{a,b}'\right)^{-1}({\mathcal Z_M'}^{a,b})\right).$$
We will omit superscript $a,b$ if there is no confusion.

In the rest of this section, we will prove the following lemma whose proof can be given by a slight modification of the proof of Proposition 2.6 in \cite{MT92}. To make this paper self contained, we will include the proof.

\bl\label{main-tech}
Let $f:\Omega\to \Omega'$ be a proper holomorphic map. Suppose there exists a point $x\in \Omega$ and a finite collection $\Omega_{\sigma_1},\ldots,\Omega_{\sigma_k}$ of characteristic subdomains of rank $1$ passing through $x$ with $\sigma_j\in Dom(f^\sharp_1)$, $j=1,\ldots,k$ such that 
\beq\label{max cond}
\bigcap_{j} X'_{f^\sharp_1(\sigma_j)}=\{f(x)\},
\eeq
where 
$$X'_{f^\sharp_1(\sigma)}=\pi'_{a,b}\left((\rho_{a,b}')^{-1}(f^\sharp_1(\sigma))\right)$$
for $(a,b)=(a_1, b_1)$.
Then $f$ has a rational extension $\widehat f:X\to X'.$
\el

We may assume that $x=0\in \Omega$. 
Let 
$\mathcal N$ be the regular locus of the closure of $f^\sharp_1(\mathcal Z_0\cap Dom(f^\sharp_1)).$
Then $\mathcal N$ is a complex manifold in $\mathcal Z_{f(0)}'$.
Since $f^\sharp_1$ is rational, $\mathcal N\cap f^\sharp_1(\mathcal Z_0\cap Dom(f^\sharp_1))$ is dense in $f^\sharp_1(\mathcal Z_0\cap Dom(f^\sharp_1))$. Therefore 
by condition \eqref{max cond}, we obtain
$$ \bigcap_{L\in \mathcal N}L=\{f(0)\}.$$
By Lemma~\ref{max rank}, there exists an integer $m>0$ and $\sigma_0\in Dom(f^\sharp_1)\cap \mathcal Z_0$ such that $A^{(m)}(\cdot;\sigma_0)$ is of maximal rank, where $A^{(m)}(\cdot;\cdot)$ is the system given in Lemma~\ref{max cond} for $h=f^\sharp_1$.
We will show that $f$ extends meromorphically to $X_{\sigma_0}\subset X$.

Let $y\in X_{\sigma_0}$. Then $\sigma_0\in \mathcal Z_y$. Consider the equation
$$ Q(u,f^\sharp_1(\sigma))=0, \quad (u, \sigma)\in \mathbb P^{N'}\times Dom(f^\sharp_1),$$
where $Q$ is a quadratic defining equation of the tautological $\mathbb P^\ell$-bundle over $\mathscr G(\ell, N').$
Then the derivatives of the above equation with respect to $\sigma\in  \mathcal Z_y\cap Dom(f^\sharp_1)$ at $\sigma_0$ depends holomorphically on a finite jet of $f^\sharp_1$ at $\sigma_0$ and tangent vectors of $\mathcal Z_y$ at $\sigma_0$.
Let
$$h_y:=f^\sharp_1\big|_{\mathcal Z_y}$$
and let
$A^{(m)}_y(\cdot, \sigma_0)$ be the affine system defined in section \ref{am} with respect to $h_y$. Note that the solution of the $A^{(m)}_y(\cdot; \sigma_0)=0$
depends meromorphically on the coefficients of the system $A^{(m)}_y(\cdot;\sigma_0)$. Since 
$\mathcal U_1(X)$ is biholomorphic to
$$\{(\sigma, T_\sigma \mathcal Z_y):(y,\sigma)\in \mathcal U_1(X)\}\subset Gr(m, T\mathcal D_1(X)),\quad m=\dim \mathcal Z_y,$$ 
$T_{\sigma_0}\mathcal Z_y$ depends holomorphically on $y$. Furthermore, $f^\sharp_1\big|_{\mathcal Z_y}$ depends meromorphically on a finite jet of $f$ at $y$. Therefore the coefficients of $A^{(m)}_y(\cdot;\sigma_0)$ depends meromorphically on $y$, implying that the solution to
$$ A^{(m)}_y(\cdot, \sigma_0)=0$$
depends meromorphically on $y\in X_{\sigma_0}$.

Note that the condition on the rank of $A^{(m)}_y(\cdot;\sigma)$ is generic. Therefore there exists a dense open set $\mathscr W_0\subset \mathcal D_1(\Omega)$ such that
$f$ extends meromorphically on $X_\sigma$ for all $\sigma\in \mathscr W_0$.

\bl 
$f$ extends rationally to $X$.
\el

\bpf
Let
$$W_1:=\bigcup_{\sigma\in \mathscr W_0}X_\sigma.$$
First we will show that $f$ extends meromorphically to $W_1$. Fix a point $\sigma_0\in \mathscr W_0$ and choose a point $y_0\in  \Omega^c\cap X_{\sigma_0}$. It is enough to show that the solution to
$$Q(\cdot, L)=0, \quad L\in f^\sharp_1(\mathscr W_0\cap \mathcal Z_y)$$ 
depends meromorphically on $y$ on an open neighborhood of $y_0$.
Choose a small ball $B\subset X$ centered at $y_0$ and choose a holomorphic section $s:B\to \mathcal U_1(X)$ such that $s(y_0)=(y_0, \sigma_0)$. After shrinking $B$ if necessary, we may assume $\rho_1\circ s(B)\subset \mathscr W_0$. Write
$$\sigma(y)=\rho_1\circ s(y),\quad y\in B.$$
Since $f^\sharp_1$ is rational, the $m$-jet of $f^\sharp_1$ at $s(y)$ depends meromorphically $y$. Hence the coefficients of the system 
$A^{(m)}_y(\cdot;\sigma)$ at $\sigma=\sigma(y)$ depends meromorphically on $y$.
It implies that the unique solution of the system 
$$Q(\cdot, L)=0, \quad L\in f^\sharp_1(\mathscr W_0\cap \mathcal Z_y)$$  
should be the meromorphic extension of $f$ on $B$.

Suppose $W_{i-1}\subset X$ and $\mathscr W_{i-1}\subset\mathcal D_r(X)$ are defined such that $f$ extends meromorphically on $W_{i-1}$ and $A^{(m)}(\cdot;\sigma)$ has maximal rank for all $\sigma\in \mathscr W_{i-1}$, where we let $W_0=\Omega$.
Define
$$W_i:=\bigcup_{\sigma\in \mathscr W_{i-1}}X_{\sigma}.$$
Then as above, $f$ extends meromorphically on $W_i$. Define
$\mathscr W_i$ to be the set of points $\sigma$ in $\mathcal Z_y, ~y\in W_i$ 
such that $A^{(m)}_y(\cdot;\sigma)$ has maximal rank at $\sigma$. Since $y\in W_i$, by definition, there exists a point in $\mathcal Z_y$ at which $A^{(m)}_y(\cdot;\sigma)$ is of maximal rank. Since the rank condition is generic, $\mathscr W_i\cap \mathcal Z_y$ is dense in $\mathcal Z_y$. In particular, 
$$\mathscr W_{i-1}\subsetneq \mathscr W_{i},$$
if $W_{i-1}$ is a proper subset of $ W_i$.

Finally we will show that there exists $i$ such that $W_i=X$. Observe that by Polysphere Theorem (\cite{M89}), for any two point $x, y\in X$, there exists a chain of $X_\sigma, ~\sigma\in \mathcal D_1(X)$ of length $\leq q$ that connects $x$ and $y$. Let $x\in \Omega, ~y\in \Omega^c$. Choose a chain $X_{\sigma_1},\ldots, X_{\sigma_a}$ connecting $x$ and $y$. 
Let 
$$ \mathscr W:=\{\sigma\in \mathcal D_r(X): A^{(m)}_z(\cdot;\sigma)\text{  is of maximal rank for some $z\in X_\sigma$.  }\}.$$
Then $\mathscr W$ is Zariski dense in $\mathcal D_1(X)$. Note that the projection map in \cite{HoN21} is a holomorphic surjection. Therefore by moving $\sigma_j$ sufficiently small, we may choose another characteristic subspaces $\widetilde\sigma_{1},\ldots,\widetilde \sigma_a$ in $\mathscr W$ whose chain connects $y$ and a point $\widetilde x\in \Omega$, i.e. $y\in W_a$. Therefore $f$ extends meromorphically on a neighborhood of $y$, which completes the proof.
\epf

\section{CR structures of boundary orbits }

The rank $r$ boundary orbit $S_r=S_r(\Omega)$ of a bounded symmetric domain $\Omega$ is a homogeneous CR manifold foliated by boundary components of rank $r$. In this section, we investigate the CR structures of $S_r=S_r(\Omega)$ and submanifolds in it for type one bounded symmetric domains $\Omega=D_{p,q}$. We refer \cite{KZ15} as a reference.

\subsection{CR structure of $S_r$}
For each $x\in S_{r}=S_{p,q,r}$, there exists a unique maximal subspace $Z_0(x)\subset V_x$ of dimension $q-r$ such that 
$$\langle Z_0(x), Z_0(x)\rangle=0.$$
We call $Z_0(x)$ the \emph{maximal null space} of $V_x$. Choose a complementary subspace $\widetilde Z\subset V_x$ of $Z_0(x)$. 
We will denote $x$ by
$$x=Z_0(x)\oplus \widetilde Z.$$
An \emph{adapted $S_{r}$-frame} or simply an $S_r$-\emph{frame} is a set of vectors
$$Z_{1},\ldots, Z_{q-r}, \widetilde Z_1,\ldots, \widetilde Z_{r}, X_{1},\ldots, X_{p-q+r}, Y_{1},\ldots, Y_{q-r}$$
in $\mathbb {C}^{p+q}$
for which the basic form $\langle~,~\rangle=\langle~, ~\rangle_{p,q}$
is given by the matrix
$$
\begin{pmatrix}
0&0& 0& I_{q-r}\\
0&-I_{r}&0&0\\
0&0&I_{p-q+r}&0\\
I_{q-r}&0&0&0\\
\end{pmatrix}.
$$
Thus we have
$$Z_0(x)=\text{span} \{Z_{1},\ldots,Z_{q-r}\},
\quad  V_x= Z_0 \oplus \text{span} \{\widetilde Z_{1},\ldots, \widetilde Z_{r}\}.$$
Denote
$$\widetilde Z:=\text{span}\{\widetilde Z_{1},\ldots,\widetilde Z_{r}\}, \quad
X:=\text{span} \{X_{1},\ldots,X_{p-q+r}\}, \quad
Y:=\text{span} \{ Y_{1},\ldots, Y_{q-r} \}.$$
Then the basic form $\langle~, ~\rangle$ defines the natural duality pairings
$$Z_0(x)\times Y\to \C,\quad \widetilde Z\times \widetilde Z\to \C, \quad X\times X\to \C.$$
 
Denote by $\mathcal B_r=\mathcal B_{p,q,r}\to S_r$ the adapted $S_r$-frame bundle
and by $\pi$ the Maurer-Cartan (connection) form on $\mathcal B_r$
satisfying the structure equation $d\pi= \pi\wedge \pi$.
Then we can write
\begin{equation}
\begin{pmatrix}
dZ_{\a}\\
d\widetilde Z_{u}\\
dX_{k}\\
dY_{\a}
\end{pmatrix}
= \pi
\begin{pmatrix}
Z_{\b}\\
\widetilde Z_{v}\\
X_{j}\\
Y_{\b}
\end{pmatrix}
=
\begin{pmatrix}
\psi_{\a}^{~\b} & \widetilde \theta_{\a}^{~v} & \theta_{\a}^{~j} & \phi_{\a}^{~\b}\\
\widetilde\sigma_{u}^{~\b} &  \widetilde\omega_{u}^{~v} & \delta_{u}^{~j} & \widetilde\theta_{u}^{~\b} \\
\sigma_{k}^{~\b} & \delta_{k}^{~v} & \omega_{k}^{~j} & \theta_{k}^{~\b}\\
\xi_{\a}^{~\b} & \widetilde\sigma_{\a}^{~v} & \sigma_{\a}^{~j} & \3\psi_{\a}^{~\b}\\
\end{pmatrix}
\begin{pmatrix}
Z_{\b}\\
\widetilde Z_{v}\\
X_{j}\\
Y_{\b}
\end{pmatrix}.
\end{equation}
In the sequel, as in \cite{KZ15}, we will identify forms on
$\mathcal B_r$ with their pullbacks to $S_r$ via local sections of the frame bundle $\mathcal B_r\to S_r$.
With that identification in mind,
the forms $\phi_{\a}^{~\b}$ give a basis
in the space of all contact forms, i.e.,
$$T^{1,0}_xS_r=\{\phi_\a^{~\b}=0,\forall \a,\b\}$$
and the upper right block forms
$$\begin{pmatrix}
\theta_{\a}^{~j} & \phi_{\a}^{~\b}\\
\delta_{u}^{~k} & \widetilde\theta_{u}^{~\b}
\end{pmatrix}
$$
form a basis
in the space of all $(1,0)$ forms on $S_r$.
We denote by $\phi$, $\theta$, $\delta$ the spaces of one forms spanned by $\{\phi_\a^{~\b},~\forall \a,\b\}$, $\{\widetilde\theta_u^{~\b},\theta_\a^{~j},~\forall \a,\b,u,j\}$ and $\{\delta_u^{~j},~\forall u, j\}$, respectively.
\medskip

There are several types of frame changes.

\bd
{\rm We call a change of frame}
\begin{enumerate}
\item[i)]change of position {\rm if}
$$
Z_\alpha'=W_\alpha^{~\beta}Z_\beta,\quad \widetilde Z'_u=W_u^{~\beta}Z_\b+W_u^{~v}\widetilde Z_v, \quad
Y_\alpha'=V_\alpha^{~\beta}Y_\beta+V_\a^{~v}\widetilde Z_v,\quad X_j'=X_j,
$$
{\rm where $W_0=(W_\alpha^{~\beta})$ and $V_0=(V_\alpha^{~\beta})$ are
$(q-r)\times (q-r)$ matrices satisfying $V_0^*W_0=I_{q-r}$, $\widetilde W=(W_u^{~v})$ is an $r\times r$ matrix satisfying $\widetilde W^*\widetilde W=I_{r}$ and $V_\a^{~\b}{W^{*}}_\b^{~\gamma}+V_\a^{~v}{W^{*}}_v^{~\gamma}=0$};

\item[ii)]change of real vectors {\rm if}
$$
 Z_\alpha'=Z_\alpha,\quad \widetilde Z'_u=\widetilde Z_u,\quad
X_j'=X_j,\quad
Y_\alpha'=Y_\alpha+H_\alpha^{~\beta}Z_\beta,
$$
{\rm where $H=(H_\alpha^{~\beta})$ is a skew hermitian matrix};

\item[iii)]dilation {\rm if}
$$
 Z_\alpha'=\lambda_{\alpha}^{-1}Z_\alpha,\quad \widetilde Z'_u=\widetilde Z_u, \quad
 Y_\alpha'=\lambda_\alpha Y_\alpha,\quad
 X_j'=X_j,
$$
{\rm where $\lambda_\alpha>0$};

\item[iv)]rotation {\rm if}
$$
Z_\alpha'=Z_\alpha,\quad \widetilde Z'_u=\widetilde Z_u,\quad
 Y_\alpha'=Y_\alpha,\quad
 X_j'=U_j^{~k}X_k,
$$
{\rm where $(U_j^{~k})$ is a unitary matrix.}
\end{enumerate}
\ed
The remaining frame change is given by
$$
 Z_\alpha'=Z_\alpha,\quad \widetilde Z'_u=\widetilde Z_u,\quad
 X_j'=X_{j} + C_j^{~\beta}Z_\beta,\quad
 Y_\alpha'=Y_\alpha+A_\alpha^{~\beta}Z_\beta+B_\alpha^{~j}X_j,
$$
such that
$$C_j^{~\alpha}+B_j^{~\alpha}=0$$
and
$$(A_\alpha^{~\beta} + \overline{A_\beta^{~\alpha}})
+B_\alpha^{~j}B_j^{~\beta}=0,$$
where
$$B_j^{~\alpha}:=\overline{B_\alpha^{~j}}.$$
The change of connection one form $\pi$ under each frame change is described in \cite{KZ15}.
\medskip

A boundary component $\Omega_\sigma$ contains $x\in S_r$ if and only if 
$$\Omega_\sigma=[Z_0(x), Z_0^*(x)]_q\cap S_r.$$
Therefore $\Omega_\sigma\subset S_r$ is a maximal integral manifold of
$$dZ_\a=0\mod Z_0\oplus \widetilde Z,\quad\forall \a$$
or equivalently, a maximal integral manifold of 
$$\phi_\a^{~\b}=\widetilde\theta_u^{~\b}=\theta_\a^{~j}=0,\quad \forall \a,\b,u,j.$$

%
%
%
%

\bl\label{R-orthogonal}
Let $A_x\subset \mathscr C_x(X)$ be the set of characteristic directions orthogonal to $T_x^{1,0}S_r$ with respect to the canonical K\"{a}hler-Einstein metric.
Then
$$T_x [Z_0(x), Z_0(x)^*]_q=\mathscr N_{A_x},$$
where $\mathscr N_A$ is the space defined in Section \ref{cb}.
\el

\bpf
Let $\{Z_\a, \widetilde Z_u, X_j,Y_\a\}$ be an $S_r$-frame at $x$. Since 
$\phi_\a^{~\b}$, $\a, \b=1,\ldots,q-r$ span the space of contact forms, after a suitable frame change, we may assume that 
$$\text{span}\{A_x\}=\textup{Hom}(Z_0(x), Y).$$
Choose a $q-r$ dimensional polyshere $P_x=\prod_{\a=1}^{q-r} \mathbb P_\a$ passing through $x$, where $\mathbb P_\a$ is a projective line such that 
$$T_x\mathbb P_\a=\textup{Hom}(Z_\a, Y_\a).$$
Then for any $q-r$ dimensional null space
$$V_0=\text{span}\{Z_\a+c_\a Y_\a,~\a=1,\ldots, q-r\},\quad c_\a\in \mathbb R,$$
and $r$ dimensional negative definite space
$$\widetilde V \subset \text{span}\{\widetilde Z_u, X_j,~u=1,\ldots, r,~j=1,\ldots, p-q+r\},$$
the point
$V_0\oplus \widetilde V$ is contained in $S_r$. 
Hence 
$\prod_{\a=1}^{q-r}\Delta_\a $ is a totally geodesic polydisc in $\Omega$ such that 
$$\prod_{\a=1}^{q-r}\partial\Delta_\a\times\Omega_\sigma\subset S_r,$$
where 
$\Delta_\a=\mathbb P_\a\cap \Omega$ and $\Omega_\sigma=[Z_0(x), Z_0(x)^*]_q\cap S_r$,
completing the proof.
\epf
Note that  by \eqref{o-decomp}, $T^{1,0}S_r$ has an orthogonal decomposition
$$T^{1,0}_xS_r=\mathscr H_{A_x}\oplus \mathscr N_{A_x}.$$
The spaces $\text{span}\{A_x\}$, $\mathscr H_{A_x}$ and $\mathscr N_{A_x}$ are parallel along $x\in \Omega_\sigma$ for each fixed boundary component $\Omega_\sigma\subset S_r$ in Harish-Chandra coordinates(cf. \cite{MT92}).

\subsection{CR structure of $\mathscr S_M^r$}
Let $M$ be a subgrassmannian of the form $M=[V_M, V_M^*]_q$ for some $q-s$ dimensional null space $V_M$ with $s<r$.
Then for any $x\in  S_r$, $x$ is contained in $\mathscr S_M=\mathscr S_M^r$ if and only if $ Z_0(x)\subset V_M,$
i.e.,
$$\mathscr S_M=\{x\in S_r: \langle Z_0(x), V_M^*\rangle=0\}.$$
Hence 
$$\text{span}\{Z_0(x):x\in \mathscr S_M\}=V_M$$
and
$\mathscr S_M\subset S_r$ is a maximal integral manifold of 
$$\langle dZ_\a, V_M^*\rangle=0,\quad \a=1,\ldots, q-r.$$
By \eqref{ZM}, $\mathscr Z_M=\mathscr Z_{M}^r$ is a Schubert variety of the form 
$$\{x\in X: \dim V_x\cap V_M\geq q-r, ~\dim V_x\cap V_M^*\geq q-r+s\}$$
when $s>0$ and of the form
$$\{x\in X: \dim V_x\cap V_M\geq q-r\}$$
when $s=0$.

\bl\label{CR on A}
$\mathscr S_M$ is a CR manifold such that
$$ T^{1,0}_x \mathscr S_M=T_x^{1,0} [Z_0(x), Z_0(x)^*]_q.
$$
Moreover, if $\Omega$ is of tube type and $P\in S_0$, then for $x\in \mathscr S_P,$  
$$\left\{v-\sqrt{-1}J(v): v\in T_x \mathscr S_P\right\}=T_x^{1,0}\mathscr Z_P= T_x^{1,0}S_r,$$
where $J$ is the complex structure of $X$.
\el
\bpf
Let $\{Z_\a, \widetilde Z_u, X_j, Y_\a\}$ be an adapted $S_r(\Omega)$-frame at $x$
so that 
$$dZ_\a=\widetilde\theta_\a^{~u}\widetilde Z_u+\theta_\a^{~j}X_j+\phi_\a^{~\b}Y_\b\mod Z_0(x).$$ 
Then on $T_x\mathscr S_M$, we obtain
\beq\label{SM-tangent}
\widetilde\theta_\a^{~u}\langle \widetilde Z_u, V_M^*\rangle+\theta_\a^{~j}\langle X_j, V_M^*\rangle+\phi_\a^{~\b}\langle Y_\b,V_M^*\rangle=0,\quad \a=1,\ldots,q-r.
\eeq
Since $V_M$ is a null space containing $Z_0(x)$, after a frame change, we may assume that
$$V_M=\text{span}\{Z_\a, \widetilde Z_u-X_u,~\a=1,\ldots, q-r,~u=1,\ldots, r-s\}.$$
Then 
$$V_M^*=V_M+\text{span}\{\widetilde Z_u, X_j,~u,j>r-s\}.$$
Hence by \eqref{SM-tangent}, we obtain
$$T_x \mathscr S_M=\{ \phi_\a^{~\b}=\widetilde\theta_\a^{~u}=\Box_\a^{~j}=0,~\a,\b=1,\ldots,q-r,~j=1,\ldots,p-q+r,~u=r-s+1,\ldots,r\},$$
where
$$\Box_\a^{~j}=\widetilde\theta_\a^{~j}+\theta_\a^{~j},\quad j=1,\ldots,r-s$$
and
$$\Box_\a^{~j}=\theta_\a^{~j},\quad j=r-s+1,\ldots,p-q+r.$$

Let $v\in T^{1,0}_x\mathscr S_M$. Then $Re(v)$ and $Im(v)$ should satisfy 
$$\phi_\a^{~\b}=\widetilde\theta_\a^{~u}=\Box_\a^{~j}=0,\quad \a,\b=1,\ldots,q-r,~j=1,\ldots,n,~u=r-s+1,\ldots,r,$$
which implies
$$\phi_\a^{~\b}=\widetilde\theta_u^{~\b}=\theta_\a^{~j}=0,\quad \forall \a,\b,u,j.$$
Therefore
$$T^{1,0}_x\mathscr S_M=T_x[Z_0(x), Z_0(x)^*]_q.$$
If $M=\{P\}$ and $\Omega$ is of tube type, then 
$s=0$ and $V_P^*=V_P$. 
Hence on a neighborhood of $x\in \mathscr S_P$, $\mathscr Z_P$ is defined by
$$ \{x\in X: \dim V_x\cap V_P=q-r\},$$
which implies
$$T_x^{1,0}\mathscr Z_M
=\{\phi_\a^{~\b}=0\}=T^{1,0}_xS_r,$$
completing the proof.
\epf

%
%
%
%
%
In the proof of Lemma~\ref{CR on A}, we showed that there exists an $S_r$-frame $\{Z_\a, \widetilde Z_u, X_j, Y_\a\}$ at $x$ such that 
$\mathscr S_M$ is a maximal integral manifold of 
$$\phi_\a^{~\b}=\widetilde\theta_\a^{~u}=\Box_\a^{~j}=0,\quad u=r-s+1,\ldots, r,~j=1,\ldots,p-q+r,$$
where
$$\Box_\a^{~j}=\widetilde\theta_\a^{~j}+\theta_\a^{~j},\quad j=1,\ldots,r-s$$
and
$$\Box_\a^{~j}=\theta_\a^{~j},\quad j=r-s+1,\ldots,p-q+r.$$
We will denote 
by $\Box_M$ the space of $1$-forms spanned by $\{\Box_\a^{~j},~\forall \a, j\}.$
If $M={P}$, we will denote $\Box_M$ by $\Box_P$ for simplicity.
Then for such frames, we obtain the following lemma.
\bl\label{L-sigma}
Let $\Omega_\sigma$ be a rank $r>s$ boundary component of $\Omega$ such that $\Omega_\sigma\subset \mathscr Z_M^r$ and let $x\in \Omega_\sigma$ be a smooth point of $\mathscr Z_M^r$.
Then $T_x^{1,0}\mathscr Z_M^r$ is a subspace of $T^{1,0}_xS_r$ defined by
$$\widetilde\theta_u^{~\b}=\theta_\a^{~j}=0,\quad u, j>r-s.$$
\el

\bpf
Since $x$ is a smooth point of $\mathscr Z_M$, $\mathscr Z_M$ is locally defined by
$$\{y\in X: \dim V_y\cap V_M= q-r, ~\dim V_y\cap V_M^*= q-r+s\}$$
when $s>0$ and
$$\{y\in X: \dim V_y\cap V_M= q-r\}$$
when $s=0$,
which completes the proof.
\epf

\section{Second fundamental form of $f$} 
\subsection{Transversality}
In this subsection we show the transversality of a holomorphic map $f$ between type one bounded symmetric domains at general Shilov boundary points. More precisely, we will show the following proposition.
\bp\label{trans}
Let $f:\Omega\to \Omega'$ be a proper holomorphic map between irreducible bounded symmetric domains of type one. Suppose there exist a point $P\in S_0(\Omega)$ and an open neighborhood $U$ of $P$ such that $f$ extends $C^2$ to $U$. Suppose further that $f$ maps $S_0(\Omega)\cap U$ to $S_0(\Omega')$.
Then for general point $x\in  S_0(\Omega)\cap U$, 
the rank of $f_*(v)$ is equal to the rank of $\Omega'$ for general $v\in T_x X$.
\ep

First we will show the following.

\bl\label{trans-disc}
Let $\Delta\subset \mathbb C$ be a unit disc and let $h:\Delta\to \mathbb B^n$, $n\geq 2$ is a proper holomorphic map, $C^2$ up to a connected open set $U$ of the boundary. Then the set $\{P\in U:h_*(\nu_P)\in T^c\mathbb \partial \mathbb B^n\}$ is a nowhere dense subset of $\partial \Delta$, where $\nu_P$ is an outward normal vector of $\Delta$ at $P$ and
$T^c\partial \mathbb B^n=T\partial \mathbb B^n\cap J(T\partial\mathbb B^n).$
\el

\bpf
Since everything in the proof is purely local, we may assume that $h$ is holomorphic on the left half plane $\mathbb H=\{\zeta=r+\sqrt{-1}t:r<0\}$ and 
$\nu_P=\dfrac{\partial}{\partial r}.$
Suppose there exists an open set $U\subset \partial \mathbb H$ such that $h$ is $C^2$ and 
$$h_*(\nu_P)\in T^c\mathbb B^n,\quad P\in U$$
or equivalently
\beq\label{tangential}
	h'(P)\cdot \overline {h(P)}=\sum_j  h_j' (P)\overline {h_j(P)}=0,\quad P\in U,
\eeq
where $h=(h_1,\ldots,h_n)$ and $h_j'=\dfrac{dh_j}{d\zeta}$.
Apply $\dfrac{\partial}{\partial t}={\rm Im}\left(\dfrac{\partial}{\partial \zeta}\right)$ to obtain
\beq\label{h''}
\sum_j h_j''(P)\overline{ h_j(P)}-\sum_j \lvert h_j'\lvert^2(P)=0,\quad P\in U.
\eeq
Write
$$h(r+it)=A_0(t)+A_1(t)r+A_2(t)r^2+o(r^2).$$
Then \eqref{tangential} and \eqref{h''} imply
$$A_1(t)\cdot \overline {A_0(t)}=0,\quad
A_2(t)\cdot\overline {A_0(t)}=\Vert A_1(t)\Vert^2,\quad t\in U.$$
Therefore for $t\in U$,
\begin{eqnarray*}
	\Vert h(r+it)\Vert^2 &=& (A_0(t)+A_1(t)r+A_2(t)r^2)\cdot \overline{(A_0(t)+A_1(t)r+A_2(t)r^2)}+o(r^2)\\
	&=&1+3\Vert A_1(t)\Vert^2 r^2+o( r^2).
\end{eqnarray*}
Since $h$ is proper, 
$$\Vert A_1(t)\Vert^2\leq 0,\quad t\in U.$$ 
Then 
$$\Vert h'(P)\Vert^2=0,\quad P\in U,$$ 
implying that $h$ is constant, contradicting the assumption that $h$ is proper.
\epf

{\it Proof of Proposition~\ref{trans}:} Choose a totally geodesic holomorphic disc $\phi:\Delta\to \Omega$ such that 
$\phi(\partial\Delta)\subset S_0(\Omega)$ and $\phi(\partial\Delta)\cap U\neq\emptyset$. Write $F:=f\circ \phi=(F_1,\ldots, F_{q'}):\Delta\to \Omega'$ in Harish-Chandra coordinates, where $F_j$ is a vector-valued function. Since $f$ maps Shilov boundary to Shilov boundary, on $\partial \Delta$, we obtain
\beq\label{diff-disc}
F_j\cdot \overline F_k=\delta_{jk},\quad j,k=1,\ldots,q'.
\eeq
By Lemma~\ref{trans-disc}, there exists $P\in \partial\Delta$ such that
$$F_j'(P)\cdot \overline F_j(P)\neq 0,\quad j=1,\ldots,q'.$$
After composing an automorphism of $\Omega'$, we may assume that 
$$F(P)=(Id_{q'};0)^t.$$
Write
$$F=(f_1;f_2)^t,$$
where $f_1$ is a ${q'}\times {q'}$ matrx-valued holomorphic function.
Then by differentiating \eqref{diff-disc} along $\partial\Delta$, we obtain
$$(f_1')(P)= (\overline{ f_1'(P)})^t,$$
i.e., $(f_1'(P))$ is a Hermitian matrix such that each diagonal entry is nonvanishing. 
This property holds up to any automorphism of $\Omega'$ that fixes $F(P)$. Hence $f_1'$ is of rank $q'$, which completes the proof.

Now let $f:\Omega\to \Omega'$ be as in the Theorem~\ref{main}. Then $f$ restricted to each boundary component is a holomorphic map into a boundary component of $\Omega'$ and extends $C^2$ up to the boundary sending Shilov boundary to Shilov boundary. Hence by Proposotion \ref{trans}, we obtain the following lemma.

From now on, we denote by $S_r(\Omega)$ and $S_r(\Omega')$, the rank $r$ boundary orbits of $\Omega$ and $\Omega'$, respectively.

\bl\label{unique bd}
Suppose there exists a nonempty open neighborhood $U$ of a Shilov boundary point of $\Omega$ such that $f$ extends $C^2$ to $U$ and $f(S_1(\Omega)\cap U)\subset S_r(\Omega')$.  
Then there exists $P\in S_0(\Omega)\cap U$ such that
$f_*(\nu)$ is of rank $r$ for general $[\nu]\in \mathscr C_P(X)$. 
\el
Since $f$ is holomorphic in $\Omega$, if $f$ satisfies the condition in Lemma~\ref{unique bd}, then $f_*(\nu)$ is of rank $r$ for general $P\in\overline{ \Omega}$ and general $[\nu]\in \mathscr C_P(X)$.
From now on to the rest of this section, we assume that $\Omega$ is of tube type and $f_*(\nu)$ is of rank $r$ for general $P\in \overline\Omega$ and general $[\nu]\in \mathscr C_P(X)$.
\medskip

%


Let  $f_1^\sharp:\mathcal D_1(X)\to \mathcal F_{a,b}(X')$ be the moduli map induced by $f$.
Since $f_*(\nu)$ is of rank $r$ for general $\nu\in \mathscr C(\Omega)$, for $\sigma\in \mathcal D_1(S_1(\Omega)),$
$f^\sharp_1(\sigma)$ is of the form
$[A, B]_{q'}$ for some $(q'-r)$-dimensional null space $A$. Since $f^\sharp_1$ is rational and $\mathcal D_1(S_1(\Omega))$ is a CR manifold that is not contained in any complex subvariety of $\mathcal D_1(X)$ of positive codimension (See \cite{KMS22}),
for any $\sigma\in Dom(f^\sharp_1)$, $f^\sharp_1(\sigma)$ is of the form
$[A, B]_{q'}$ for some $(q'-r)$-dimensional subspace $A$. That is, $f^\sharp_1$ is a map to $\mathcal F_{q'-r, b}(X')$ for some $b$.

Let
$$Pr:\mathcal F_{q'-r, b}(X')\to Gr(q'-r, \mathbb C^{p'+q'})$$
be a projection defined by $Pr((A, B))=A$.
For a general point $P\in S_0(\Omega)$, define
$$\mathcal L_P:=\{(A, A^*)\in \mathcal D_r(X'):A\in Pr(\overline{f_1^\sharp(\mathcal S_P)\cap Dom(f_1^\sharp)})\}$$
and
\beq\label{scr LP}
\mathscr L_P:=\pi_r'\left(\left(\rho_r'\right)^{-1}\left(\mathcal L_P\right)\right),
\eeq
where $\pi'_r:\mathcal U_r(X')\to X'$, $\rho_r':\mathcal U_r(X')\to \mathcal D_r(X')$ is the canonical double fibration of the universal family of characteristic spaces of rank $r$ over $X'$ defined in Section \ref{am}. 
Since $f$ is a CR map,
for each general boundary component $\Omega_\sigma\subset S_1(\Omega)$, there exists a unique boundary component of $\Omega'$ with rank $r$ that contains $f(\Omega_\sigma)$.  This boundary component is given by
$$[Z_0(f(x)), Z_0(f(x))^*]_{q'}\cap S_r(\Omega'),\quad x\in \Omega_\sigma,$$
where we denote by $Z_0(y)$ the maximal null space of $V_y$ for $y\in \partial \Omega'$. Recall that $Z_0(f(x))$ is constant on each boundary component.  
Therefore
$$Z_0(f(\mathscr S_P\cap U)))\subset Pr(\overline{f_1^\sharp(\mathcal S_P)\cap Dom(f_1^\sharp)})$$
and
hence
\beq\label{def-LP}
f(\mathscr S_P\cap U)\subset \mathscr L_P.
\eeq
Since each point in the boundary orbit of rank $r$ has a unique maximal null space of codimension $r$, we can define a smooth map
$\Pi:S_r(\Omega')\to Gr(q'-r, \mathbb C^{p'+q'})$ by $\Pi(y)=Z_0(y).$ Then for any general $x\in \mathscr S_P\cap U$, there exists an open neighborhood of $U'\subset S_r(\Omega')$ of $f(x)$ such that 
$$ \mathscr L_P\cap U'=\Pi^{-1}\left(Z_0(f(\mathscr S_P^1\cap U))\right)\cap U'.$$
For a general point $x\in \mathscr S_P\cap U$, we may assume that $\Pi\circ f$ is a surjection at $x$ onto its image
and hence
$$T_{f(x)}\mathscr L_P=f_*(T_x\mathscr S_P^1(\Omega))+ker(\Pi_*).$$
By taking the complexification, we obtain
\beq\label{MP-tangent}
\{v-\sqrt{-1}Jv:v\in T_{f(x)}\mathscr L_P\}
=f_*(T^{1,0}_x S_1)+T_{f(x)}[Z_0(f(x)), Z_0(f(x))^*]_{q'}.
\eeq
\medskip

Let $P\in S_0(\Omega)\cap \overline\Omega_\sigma$ for some rank one boundary component $\Omega_\sigma$ and let $x=P+\lambda \nu_\sigma\in \Omega_\sigma\subset S_1(\Omega)$ be general points in $U$, where $\nu_\sigma\in T_PX$ is an outward normal vector of $\partial\Omega_\sigma$ in $X_\sigma$ at $P$ and 
$\lambda<0$ is a sufficiently small real number. Then
$$f(P+\lambda \nu_\sigma)=f(P)+\lambda  f_*(\nu_\sigma)+O(\lambda^2)\in \Omega'_{\sigma'},$$
where $\Omega_{\sigma'}'$ is the rank $r$ boundary componentof $\Omega'$ that contains $f(\Omega_\sigma\cap U)$.
Since $f_*(\nu_\sigma)$ is of rank $r$ and tangent to $f(X_\sigma)$ at $f(P)$, we may assume that 
$f(P+\lambda\nu_\sigma)$ and $f(P)+\lambda f_*(\nu_\sigma)$ are contained in the same boundary component $\Omega'_{\sigma'}$ for all sufficiently small $\lambda<0$, i.e.,
$$ 
Z_0(f(x))=Z_0(f(P)+\lambda f_*(\nu_\sigma)).
$$
Since $T_y\mathscr L_{P}$ is parallel along $y\in [Z_0(f(x)), Z_0(f(x))^*]_{q'}$ in Harish-Chandra coordinates, 
it implies 
$$T_{f(x)} \mathscr L_P=T_{f(P)+\lambda f_*(\nu_\sigma)}\mathscr L_P$$
under parallel translation.

Now consider a one parameter family of rank one boundary components $\Omega_{\sigma(t)}$, $t\in (-\epsilon,~\epsilon)$ such that
$$T_PX_{\sigma(t)}=(a+tb)^*\otimes (a+tb).$$
That is, after a suitable frame change, we can choose an $S_1(\Omega)$-frame $\{e_1,\ldots,e_q, e_1^*,\ldots,e_q^*\}$ of $\mathbb C^{2q}$ with
$$\langle e_i,e_j\rangle=\langle e_i^*,e_j^*\rangle=\langle e_i,e_j^*\rangle-\delta_{i,j}=0$$
such that 
$$V_P=\text{span}\{e_1,\ldots,e_q\}$$
and
$$X_{\sigma(t)}=[V_{\sigma(t)}, V^*_{\sigma(t)}]_q,$$
where
$$V_{\sigma(t)}=\text{span}\{e_1+\sqrt{-1}t e_q, e_2,\ldots,e_{q-1}\},\quad
V_{\sigma(t)}^*=\text{span}\{e_q,e_q^*+\sqrt{-1}te_{1}^*\}+V_{\sigma(t)}.$$
Choose two curves
\beq\label{xt}
x(t):=f(P+\lambda\nu_{\sigma(t)}),\quad t\in (-\epsilon, ~\epsilon)
\eeq
and 
$$
	y(t):=f(P)+\lambda f_*(\nu_{\sigma(t)})
	=f(P)+\lambda\left( f_*(\nu_a)+t f_*(v_{a,b})+ t^2 f_*(\nu_b)\right),\quad t\in (-\epsilon, ~\epsilon),
$$
where $\nu_a,\nu_b $ and $v_{a,b}$ are vectors in $T_PX$ such that the vector field 
$$\nu_{\sigma(t)}:=\nu_a+tv_{a,b}+t^2\nu_b,\quad t\in (-\epsilon,~\epsilon)$$
satisfies
$$T_PX_{\sigma(t)}=\mathbb C \nu_{\sigma(t)}. $$
After shrinking $(-\epsilon,~\epsilon)$, we may assume that for a sufficiently small fixed $\lambda<0$, 
$$Z_0(x(t))=Z_0(y(t)),\quad t\in (-\epsilon,~\epsilon).$$

Set
$$\mathscr L_{a,b}:=\bigcup_{t\in (-\epsilon,~\epsilon)} [Z_0(y(t)), Z_0(y(t))^*]_{q'}.$$
Since $f(P)$ is a Shilov boundary point, $Z_0(y(t))$ is a subspace of $V_{f(P)}$, i.e.,
$$Z_0(y(t))=V_{y(t)}\cap V_{f(P)}.$$
Therefore $\mathscr L_{a,b}$ is a submanifold of a Schubert variety 
$$\mathcal W_P:=\{y\in Gr(q', p'):\dim V_y\cap V_{f(P)}\geq q'-r\}.$$
Moreover, $y(0)$ is a smooth point of $\mathcal W_P$. 
On the other hand, since $y(t)$ is a curve of degree two, 
$\{y(t):t\in (-\epsilon,~\epsilon)\}$ is contained in 
$$y(0)+\text{span}\{ \dot{y}(0),\ddot{y}(0)\}\subset y(0)+T_{y(0)}\mathcal W_P,$$
where
$$\dot{y}(0)=\frac{d y}{ d t}(0)=f_*(v_{a,b}),\quad \ddot{y}(0)=\frac{d^2 y}{ d t^2}(0)=2f_*(\nu_b).$$
Here we regard $y(0)+T_{y(0)}\mathcal W_P$ as a linear subset of $X'$ passing through $y(0)$.
Therefore the curve $Z_0(y(t))$, $t\in (-\epsilon,~\epsilon)$ is of degree two in $t$ and contained in a linear subspace
$$Z_0(y(0))+\text{span}\{\Pi_*(\dot{y}(0)), \Pi_*(\ddot{y}(0))\}.$$
Here we regard $Z_0(y(0))+\text{span}\{\Pi_*(\dot{y}(0)), \Pi_*(\ddot{y}(0))\}$ as a linear subset of $Gr(q'-r, V_{f(P)})$ passing through $Z_0(y(0))$.
Since 
$$\Pi_*(\dot{y}(0))=\frac{dZ_0(y(t))}{dt}(0),\quad \Pi_*(\ddot{y}(0))=\frac{d^2 Z_0(y(t))}{dt^2}(0),$$
we obtain 
$$ Z_0(x(t))=Z_0(y(t))\in Z_0(f(x))+\text{span}\left \{(Z_0)_*(\dot {x}(0)), \frac{d^2 Z_0\circ x(t)}{dt^2}(0) \right\},\quad t\in (-\epsilon,~\epsilon).$$
Since $a, b$ are arbitrary, we obtain
$$ Z_0(f(\mathscr S_P\cap U))\subset Z_0(f(x))+T_{Z_0(f(x))}Z_0(f(\mathscr S_P))+\mathbb{FF}^{(2)}_{Z_0(f(x))}Z_0(f(\mathscr S_P)),$$
where 
$\mathbb{FF}^{(2)}_{Z_0(f(x))}Z_0(f(\mathscr S_P))$ is the span of the second fundamental form of $Z_0(f(\mathscr S_P\cap U))$ at $Z_0(f(x))$ with respect to the flat(Euclidean) connection of $T_{Z_0(f(x))}Gr(q'-r, V_{f(P)})\subset Gr(q'-r, V_{f(P)})$ in a big Schubert cell.

%
%
%
%
%
%
%
%
%
%
%
%

\subsection{Second fundamental form of $f:S_1(\Omega)\cap U\to S_r(\Omega')$} We will use capital Greek letters 
$\Phi_\a^{~\b}, \widetilde\Theta_U^{~\b},\Theta_\a^{~J}, \Delta_U^{~J},$ etc. for connection one forms on $\mathcal B_r(\Omega')=\mathcal B_{p',q',r}$ pulled back to $S_r(\Omega')=S_{p',q',r}$.
Since $f$ is a CR map, $f$ satisfies 
$$f^*(\Phi_\a^{~\b})=0\mod\phi,$$
$$f^*(\Theta_\a^{~J})=f^*(\widetilde\Theta_U^{~\b})=0\mod \phi, \theta,$$
\beq\label{del-in}
f^*(\Delta_U^{~J})=0\mod\phi, \theta, \delta.
\eeq
Since everything in this section is purely local and $f$ restricted to $S_1(\Omega)$ is a local embedding on an open set, we may omit $ U$ and identify $T_x^{1,0}S_1(\Omega)$ with $f_*(T_{x}^{1,0}S_1(\Omega))$. Then the pull back of one forms via $f$ is the restriction of one forms to $T_{f(x)}^{1,0}f(S_1(\Omega))$. 
In what follows, we will omit $f^*$ if there is no confusion.
\medskip

For $x\in S_1(\Omega)$, $f_*(T_x^{1,0}S_1(\Omega))$ is a subspace in ${\rm Hom}(V_{f(x)}, \mathbb C^{p'+q'}/V_{f(x)})$. We define subspaces $K_x$ and $R_x$ in $V_{f(x)}$ and $\mathbb C^{p'+q'}/V_{f(x)}$, respectively by
$$
K_{x}:=\bigcap\left\{{ {\rm Ker}}\left(proj(v)\right): v\in f_*(T_{x}^{1,0}S_1(\Omega))\right\},
\quad
R_{x}:={{\rm Span_\mathbb C}}\left\{ {\rm Im}\left(proj(v)\right): v\in f_*(T_{x}^{1,0}S_1(\Omega))\right\}.$$
where $proj$ is the projection to the orthogonal complement of $T_{f(x)}[Z_0(f(x)),Z_0(f(x))^*]_{q'}\subset T_{f(x)}X'$ with respect to the canonical K\"{a}hler-Einstein metric.
Then
$$Gr_{x}:=\left\{A\in \textup{Hom}\left(V_{f(x)}, R_{x}\right):{\rm Ker}(A)\supset K_{x}\right\}\cap T_{f(x)}^{1,0}S_r(\Omega')
+T_{f(x)}[Z_0(f(x)), Z_0(f(x))^*]_{q'}$$
is a linear subspace in $T_{f(x)}^{1,0}S_r(\Omega')$ that contains $f_*( T^{1,0}_{x}S_1(\Omega)).$

Let $P\in S_0(\Omega)$ be a point such that $x\in \mathscr S_P$. Then $f(P)$ is a point in $S_0(\Omega')$ such that
\beq\label{SP-in}
f(\mathscr S_P)\subset \mathscr S_{f(P)}'.
\eeq
After a rotation and position change, we may assume 
$$V_{f(P)}=Z_0(f(x))+\text{span}_{\mathbb C}\{  V_U:=\widetilde Z_U-X_U, ~ U=1,\ldots, r\}$$
and
$$T_{f(x)}\mathscr Z_{f(P)}'=\{\Phi_\a^{~\b}=\Theta_\a^{~J}=0,~\a,\b=1,\ldots,q'-r,~J>r\}.$$
Since $\Omega$ is of tube type and therefore
$$f_*(T^{1,0}_xS_1(\Omega))=f_*(T_x\mathscr Z_P)\subset T_{f(x)}\mathscr Z_{f(P)}',$$
we obtain
$$\Theta_\a^{~J}=0\mod \phi, \quad J>r.$$
Hence 
$$Gr_x\subset \{\Phi_\a^{~\b}=\Theta_\a^{~J}=0,~J>r\},\quad \forall \a, \b, U.$$
After rotation and position change, we may assume 
$$
Gr_x=\{\Phi_\a^{~\b}=\Theta_\a^{~J}=\widetilde\Theta_U^{~\b}=0,\quad \a,\b=1,\ldots,q'-r,~J>J_1,~U>U_1\}
$$
for some integers $U_1, J_1\leq r$. Then on $f_*(T_{x}^{1,0} S_1(\Omega))$, it holds that
$$
\Theta_\a^{~J}=\widetilde\Theta_U^{~\b}=0\mod\phi,\quad J>J_1,~U>U_1.
$$
On the other hand, by \eqref{SP-in}, we obtain
$$\widetilde\Theta_\a^{~U}+\Theta_\a^{~U}=0\mod\phi, \Box_P.$$
Hence if $U>U_1$, then
$$\Theta_\a^{~U}=0\mod\phi, \Box_P.$$
Since $f$ is a CR map, it implies 
$$\Theta_\a^{~J}=0\mod\phi,\quad J>U_1.$$
Similarly, we obtain
$$\widetilde\Theta_\a^{~U}=0\mod\phi,\quad U>J_1.$$
Therefore, 
$$J_1=U_1$$
and 
\beq\label{emb1}
\Theta_\a^{~J}=\widetilde\Theta_U^{~\b}=0\mod\phi,\quad U, J>J_1.
\eeq
Hence for general $x\in S_1(\Omega)$, there exists a reduction of $S_r(\Omega')$-frame such that
$$Gr_x=\{\Phi_\a^{~\b}=\Theta_\a^{~J}=\widetilde\Theta_U^{~\b}=0,\quad \a,\b=1,\ldots,q'-r,~J, U>J_1\}.$$
Moreover, by definition of $K_x$ and $R_x$,
for each fixed $j\leq J_1$ and $u\leq J_1$, there exist $\a$ and $\b$ such that $\Theta_\a^{~j}$ and $\widetilde \Theta_u^{~\b}$ modulo $\phi$ are nonvanishing. 
From now on, we let the small indices $u,v$, $j,k$ run from $1$ to $J_1$ and the capital indices $U, V$, $J,K$ from $J_1+1$ unless specified otherwise. 
\medskip

By differentiating \eqref{emb1} using the structure equation, we obtain
\beq\label{first prol}
\widetilde\Theta_\a^{~v}\Delta_v^{~J}+\Theta_\a^{~k}\Omega_k^{~J}=
\widetilde\Omega_U^{~v}\widetilde\Theta_v^{~\b}+\Delta_U^{~k}\Theta_k^{~\b}=0\mod \phi, \theta\wedge\overline\theta.
\eeq
We will show 
\beq\label{prol-1}
\Delta_v^{~J}=\Delta_U^{~k}=0\mod \theta, \phi
\eeq
and
\beq\label{prol-2}
\Omega_k^{~J}=\widetilde\Omega_U^{~v}=0\mod\phi, \theta, \overline\theta.
\eeq
Fix $J>J_1$.
By induction on $\a$, we can choose a position change which still satisfies \eqref{emb1} and a sequence of positive integers $v_1<v_2<\ldots<v_{q'-r}$ such that
$\widetilde \Theta_\a^{~v}, ~v_{\a-1}<v\leq v_\a$ modulo $ \phi$ is linearly independent and $\widetilde\Theta_\a^{~v}$ modulo $\phi$ vanishes for $v>v_\a$, where we let $v_0=0$. 
By condition on $J_1$, we obtain $v_{q'-r}=J_1.$
Hence by using \eqref{first prol} inductively on $\a$, we obtain
$$\sum_{v_{\a-1}<v\leq v_\a}\widetilde\Theta_\a^{~v}\Delta_v^{~J}=0\mod\phi, \theta,$$
which implies by Cartan lemma and \eqref{del-in},
$$\Delta_v^{~J}=0\mod\phi,\theta.$$
The same argument is valid for other cases.

Choose a curve $x:(-\epsilon,~\epsilon)\to \mathscr S_{f(P)}'$ of the form
\eqref{xt} such that $x(0)=f(x).$ 
We may assume that
$$\dot{x}(0)\not \in T_{f(x)}[Z_0(f(x)), Z_0(f(x))^*]_{q'}.$$ 
Write
\beq\label{quad-t}
Z_\a(x(t))=t\sum_{U=1}^r  C^{~U}_\a V_U+t^2 \sum_{U=1}^r  D_\a^{~U} V_U+O(t^3)\mod Z_0(f(x_0)),\quad \a=1,\ldots,q'-r,
\eeq
where
$$V_U=\widetilde Z_U-X_U, ~U=1,\ldots,r$$
are vectors that, together with $Z_0(f(x))$,
span $V_{f(P)}$.
Write 
$$Z_\a(x(s))=Z_\a((x(t))+(s-t)V_\a(t)+O((s-t)^2),\quad s, t\in (-\epsilon,~\epsilon)$$
for some vector field $V_\a(t)$ of the form
$$V_\a(t)=\sum_{U=1}^r C_\a^{~U}(t)V_U\mod Z_0(f(x)),\quad t\in (-\epsilon,~\epsilon).$$
Then in view of \eqref{quad-t}, 
$$C_\a^{~U}(t)=C_\a^{~U}$$
and
\beq\label{dV}
\sum_{U=1}^r  C_\a^{~U}V_U(t)=\sum_{U=1}^r C_\a^{~U}V_U+2t\sum_{U=1}^r  D_\a^{~U}V_U+O(t^2).
\eeq
By \eqref{emb1}, 
$$C_\a^{~U}=0,\quad U>J_1.$$
By differentiating \eqref{dV} with respect to $t$, 
\beq\label{range-D}
 C_\a^{~u}\dot{V_u}(0)=2\sum_{U=1}^r D_\a^{~U}V_U.
\eeq

On the other hand, since
$$dZ_\a=\frac12(\widetilde\Theta_\a^{~u}-\Theta_\a^{~u})V_u\mod\phi,\Box_{P},$$
$\text{span}\{\dot{V}_u(0),~u=1,\ldots,J_1\}$ is obtained from
\beq\label{dv}
dV_U=d\widetilde Z_U-dX_U=\sum_{W=1}^r (\widetilde \Omega_U^{~W}-\overline\Delta_U^{~W})\widetilde Z_W+\sum_{J=1}^{p'-q'+r} (\Delta_U^{~J}-\Omega_U^{~J})X_J\mod Z_0, Y.
\eeq
Define
$$\widehat V_U=\widetilde Z_U+X_U,\quad U=1,\ldots, r$$
so that
$$\langle V_U, \widehat V_W\rangle=-2\delta_U^{~W},\quad U, W=1,\ldots,r$$
and
$$\langle \widehat V_U,  X_J\rangle=0,\quad U=1,\ldots, r,~ J>r.$$
By substituting
$$\widetilde Z_W=\frac12(V_W+\widehat V_W),\quad X_W=\frac12(V_W-\widehat V_W)$$
to \eqref{dv}
we obtain
$$2dV_u
=\sum_{W=1}^r(\widetilde\Omega_u^{~W}+\Omega_u^{~W}-\Delta_u^{~W}-\overline\Delta_u^{~W})\widehat V_W+
\sum_{r>J}(\Delta_u^{~J}-\Omega_u^{~J})X_J \mod V_{f(P)},$$
which implies
\beq\label{mod-box}
(\widetilde\Omega_u^{~W}+\Omega_u^{~W}-\Delta_u^{~W}-\overline\Delta_u^{~W})
=(\Delta_u^{~J}-\Omega_u^{~J})=0,\mod\phi, \Box_P, \quad 1\leq W\leq r,~J>r
\eeq
and
$$dV_u=\sum_{W=1}^r(\widetilde\Omega_u^{~W}-\overline\Delta_u^{~W})V_W\mod \phi, \Box_P.$$

Let 
\beq\label{def-LP}
L_P:=\text{span}\{Z_0(f(y)):y\in \mathscr S_P\}.
\eeq 
Then $L_P$ is the smallest subspace of $V_{f(P)}$ such that 
\beq\label{LP-in}
\mathscr L_P\subset\pi_r'\left((\rho_r')^{-1}(\{(A, B)\in \mathcal D_r(X'): A\subset L_P, L_P^*\subset B\}\right),
\eeq
where $\mathscr L_P$ is defined in \eqref{scr LP}.
By nondegeneracy of $\widetilde\Theta_u^{~\b}, u=1,\ldots, J_1$ and \eqref{quad-t}, $L_P$ contains $\text{span}\{V_u, ~u=1,\ldots ,J_1\}$. Therefore we may assume 
$$L_P=Z_0(f(x))+\text{span}\{V_U, ~U=1,\ldots ,J_P\}$$
and
$$L_P^*=L_P+\text{span}\{\widetilde Z_U,X_J, ~U,J>J_P\}$$
for some $J_1\leq J_P\leq r.$ Since $Z_0(f(y))$ is a subspace in a fixed null space $V_{f(P)}$, 
$$L_P^*:=\bigcap_{y\in \mathscr S_P}Z_0^*(f(y)).$$
Define
$$M_P:=[L_P, L_P^*]_{q'}=\bigcap_{y\in \mathscr S_P}[Z_0(f(y)), Z_0(f(y))^*]_{q'}.$$
By definition, for any $x\in \mathscr S_P$, 
$[Z_0(f(x)), Z_0(f(X))^*]_{q'}$ contains $M_P$ and \eqref{LP-in} becomes
$$\mathscr L_P\subset\pi_r'\left((\rho_r')^{-1}_r(\mathcal Z'_{M_P})\right)=\mathscr Z_{M_P}'.$$

If $M_P$ is a point, then 
$$   \bigcap_{\sigma\in \mathcal S_P\cap Dom(f_1^\sharp)} X_{f_1^\sharp(\sigma)}'=\{f(P)\}.$$
Hence by applying Lemma~\ref{main-tech}, we can show that $f$ has a rational extension.
Now assume that $M_P$ is positive dimensional. Since $L_P$ is a a null space, $M_P$ is a nontrivial characteristic subspace of $X'$ passing through $f(P)$ such that 
$M_P\cap S_s(\Omega')$ is a boundary component of $\Omega'$ for some $0<s<r$
and
$$T_{f(P)}M_P=\bigcap_{y\in \mathscr S_P}T_{f(P)}[Z_0(f(y)), Z_0(f(y))^*]_{q'}.$$ 
By Lemma~\ref{R-orthogonal},
$$\bigcap_{y\in \mathscr S_P}T_{f(P)}[Z_0(f(y)), Z_0(f(y))^*]_{q'}=\bigcap_{y\in \mathscr S_P}\mathscr N_{A_y},$$
where 
$$A_y=\textup{Hom}(Z_0(f(y)), Y_y)$$
for some suitable $S_r(\Omega')$-frame at $f(y)$. 
Here we regard $\mathscr N_{A_y}$ as a subspace in $T_{f(P)}X'$ by parallel translation in Harish-Chandra coordinates.
Let
$\widehat L_P$ be a subspace in $\mathbb C^{p'+q'}$ such that 
$$\bigcap_{y\in \mathscr S_P}\mathscr N_{A_y}=\mathscr N_A$$
for 
$$A=\textup{Hom}(L_P, \widehat L_P/L_P).$$
That is, 
$$\widehat L_P=\text{span}\{Y_y,~y\in \mathscr S_P\}\mod L_P.$$
Since $Y_y$ is a dual of $Z_0(y)$ under the basic form $\langle~, ~\rangle$, we obtain
$$\widehat L_P=\text{span}\{\widehat V_U, ~U=1,\ldots,J_P\}\mod Y_{f(x)}+L_P.$$
and
$$T_{f(P)}M_P=\sum_{U, J>J_P}\Delta_U^{~J}X_J.$$
Therefore
$$\langle dV_u, \widetilde Z_W\rangle=\langle d\widehat V_u, X_J\rangle=0\mod \phi,\Box_P,\quad W, J>J_P,$$
which implies
$$\Delta_U^{~j}=\Delta_u^{~J}=0\mod\phi,\quad U,J>J_P$$
and
$$\widetilde\Omega_u^{~W}=\Omega_u^{~J}=0\mod\phi,\Box_P,\quad W, J>J_P.$$
Furthermore, since
$$\langle V_u, \widetilde Z_W\rangle=\langle V_u, X_J\rangle=0,\quad W, J>J_P$$
on $\mathscr L_P$,
by complexification and \eqref{MP-tangent}, we obtain
$$\langle d V_u, \widetilde Z_W\rangle=\langle d V_u, X_J\rangle=0,\quad W, J>J_P$$
on $T^{1,0}_xS_1(\Omega)$,
i.e.,
$$\langle dV_u,\widetilde Z_W\rangle=\langle d\widehat V_u, X_J\rangle=0\mod \phi,\bar\theta,\quad W, J>J_P,$$
which implies
$$
\widetilde\Omega_u^{~W}=\Omega_u^{~J}=0\mod\phi,\bar\theta,\quad W, J>J_P.
$$
Since $\Box_P$ and $\bar\theta$ are linearly independent, together with \eqref{mod-box}, we obtain
\beq\label{second-mod}
\widetilde\Omega_u^{~W}=\Omega_u^{~J}=0\mod\phi,\quad W, J>J_P.
\eeq
Moreover, since $L_P$ is the smallest subspace that satisfies \eqref{LP-in}, by the property of degree two curves, we obtain that
for each $W\leq J_P$, there exists $u$ such that 
$$\widetilde \Omega_u^{~W}-\overline\Delta_u^{~W}\neq 0\mod\phi, \Box_P,$$
implying that
$$\widetilde\Omega_u^{~W}\neq 0\mod\phi, \bar\theta.$$
This condition depends only on $\widetilde\Omega_u^{~W}$ and independent of the choice of $P$. Hence we can choose an integer $J_2=J_P$ and a further reduction of $S_r(\Omega')$-frame such that 
$$V_U=\widetilde Z_U\mod X,\quad U=J_1+1,\ldots,J_2$$
with the nondegeneracy condition on $\widetilde\Omega_u^{~W}$ modulo $\phi, \bar\theta$ stated above.
That is, if we let
$$\widetilde \Omega_u^{~W}=h^{~W,\a}_{u}\theta_\a\mod\phi,\bar\theta$$
and
$$h^{~a}_u:=\sum_{W=J_1+1,\ldots,J_2}h^{~W, a}_uV_W,$$
then
$$\text{span}\{V_u, ~u=1,\ldots, J_1\}+\text{span}\{h^{~a}_u: a=1,\ldots,q-1,~, u=1,\ldots,J_1\}=L_P\mod Z,$$
where $\theta_a,~a=1,\ldots, q-1$ are $(1,0)$-forms of $S_1(\Omega)$ corresponding to $\Theta_\a^{~J}$.

Suppose that after rotation of $X_W,~J_1<W \leq J_2,$ there exists $W_0\leq J_2$ such that 
$$\Omega_u^{~W_0}-\Delta_u^{~W_0}=0\mod\phi,\Box_P\quad \forall u.$$
Then
by \eqref{mod-box}, 
$$\widetilde\Omega_u^{~W_0}=0\mod\phi,\bar\theta,\quad \forall u,$$
which
contradictis the assumption on $J_2$.
Therefore for each $W=J_1+1,\ldots, J_2$, there exists $u$ such that 
$$\Omega_u^{~W}-\Delta_u^{~W}\neq 0\mod \phi, \Box_P,$$
which implies
$$\Omega_u^{~W}\neq0\mod\phi, \theta.$$
Thus the choice of $\widehat L_P$ is independent of $P$.
\medskip

Summing up, we have a reduction of $S_r(\Omega')$-frame such that 
$$\widetilde\Theta_U^{~\b}=\Theta_\a^{~J}=0\mod\phi,\quad U, J>J_1,$$
$$\widetilde\Omega_u^{~W}=\Delta_U^{~j}=\Delta_u^{~J}=\Omega_k^{~J}=0\mod\phi,\quad U,W, J>J_2$$
and
$$
T_{f(P)}M_P=\sum_{U, J>J_2}\Delta_U^{~J}X_J
$$
for all $P\in S_0(\Omega)$ such that $x\in \mathscr S_P$.
Here we regard $T_{f(P)}M_P$ as a subspace in $T_{f(x)}^{1,0}S_r(\Omega')$ by parallel translation in Harish-Chandra coordinates.
By using the reduction of $S_r(\Omega')$-frame, we will prove the following main technical lemma of the paper.

\bl\label{main-tech}
Let $\Omega$ be of tube type and let $P\in S_0(\Omega)$ be a general point. Suppose there exists a nontrivial subgrassmannian $M$ such that 
$$
f^\sharp_1(\mathcal S_P\cap Dom(f_1^\sharp))\subset {\mathcal S_{M}'}.
$$
Then there exists a unique maximal characteristic subspace $M_{P}\subset X'$ of the form
$$M_P=[L_P, L_P^*]_{q'}$$
for some null space $L_P$ such that 
\beq\label{char-in}
f(\mathscr S_P)\subset {\mathscr S_{M_P}'}.
\eeq
Furthermore, $M_{P}$ is parallel with $M_{\widetilde P}$ for general $P,\widetilde P\in S_0(\Omega)$.
\el

\bpf
It is enough to show that $M_P$ is parallel with each other for general $P\in S_0(\Omega)$. 
Let $x\in S_1(\Omega)$ be a general point.
Then for all $P\in S_0(\Omega)$ such that $x\in \mathscr S_{ P}^1$ or equivalently,
$$P\in [Z_0(x), Z_0(x)^*]_{q}\cap S_0(\Omega),$$
$M_P$ is parallel with each other.
Since any two points in $S_0(\Omega)$ are connected by chain of $ [Z_0(x), Z_0(x)^*]_q\cap S_0(\Omega), ~x\in S_1(\Omega)$, we can complete the proof.  
\epf


\bc\label{image}
Let $\Omega_1'$ and $\Omega_2'$ be totally geodesic subdomains of $\Omega'$ such that $\Omega_1'\times\Omega_2'$ is a totally geodesic subspace of $\Omega'$ of maximal rank passing through $f(0)$ and 
$T_{f(0)}\Omega_2'$ is parallel with $T_{f(p)}M_p$ for $P\in S_0(\Omega)$.
Then
\beq\label{img}
f(\Omega)\subset \Omega_1'\times \Omega_2'.
\eeq
Moreover, if we decompose $f=f_1\times f_2:\Omega\to \Omega_1'\times\Omega_2'$, then
$f_1$ is a proper rational map.
\ec

\bpf
Since 
$$S_0(\Omega_1'\times\Omega_2')=S_0(\Omega_1')\times S_0(\Omega_2'),$$
to show \eqref{img}, 
it is enough to show that
$$f(S_0(\Omega))\subset S_0(\Omega_1')\times S_0(\Omega_2').$$
Let $P\in S_0(\Omega)$ be a general point and let $M_P$ be the maximal characteristic subspace in Lemma~\ref{main-tech}.
Since 
$T_{f(P)}M_P$ is parallel with $T_{f(0)}\Omega_2'$,  
$M_P$ is of the form
$$M_P=\{A_P\}\times M,$$
where $M$ is the compact dual of $\Omega_2'$ and $A_P$ is a point in $S_0(\Omega_1')$.
Since $P$ is arbitrary, we obtain
$$f(P)\in S_0(\Omega_1')\times S_0(\Omega_2'),\quad \forall P\in S_0(\Omega).$$

Now it is enough to show that $f_1$ is proper and the characteristic subspace
$M_P$ in Lemma~\ref{main-tech} for $f_1$ is a point for some $P\in S_0(\Omega)$.
Suppose $f_1$ is proper. 
In the proof of Corollary~\ref{image},
$M_P$ for $f$ is of the form $\{A_P\}\times M$ for general $P\in S_0(\Omega)$.
Hence by definition, $M_P$ in Lemma~\ref{main-tech} for $f_1$ is $\{A_P\}$, implying that $f_1$ is rational.
To show that $f_1$ is proper, it is enough to show that $f_1(S_{q-1}(\Omega))\subset \partial\Omega_1'.$ Suppose otherwise. Since $f$ is proper and
$$\partial(\Omega_1'\times\Omega_2')=(\partial\Omega'_1\times\Omega_2')\cup (\Omega_1'\times\partial\Omega'_2)\cup(\partial\Omega_1'\times\partial\Omega_2'),$$ 
we obtain 
$$f_2(S_{q-1}(\Omega))\subset \partial \Omega_2',$$
i.e., $f_2:\Omega\to\Omega'_2$ is a proper holomorphic map.
Now suppose $f_2$ is proper.  
Then by Corollary~\ref{image}, there exists a further decomposition $f_2=g\times h:\Omega\to \Omega_3'\times \Omega_4'\subset \Omega_2'$ and $M_P$ in Lemma~\ref{main-tech} for $f_2$ is of the form 
$\{B_P\}\times N$ for some $B_P\in S_0(\Omega_3')$. Since $\Omega_1'\times \Omega_2'$ is totally geodesic, this implies that $M_P$ for $f$ should be of the form $\{A_P\times B_P\}\times N$, contradicting the uniqueness of $M_P=\{A_P\}\times M$ for $f$.
\epf

\section{Proof of Theorem~\ref{main} and Corollary~\ref{cor}}

First assume that $\Omega$ is of tube type. Then Corollary~\ref{image} will complete the proof.
Now assume that $\Omega$ is of non-tube type. Then the following lemma will complete the proof.

\bl
Let $P\in S_0(\Omega)$ and let 
$M_P$
be a characteristic subspace passing through $f(P)$ such that
\beq\label{EDS}
f(\mathscr S_P^1)\subset \mathscr S_{M_P}'.
\eeq
Then $M_P$ is a point for general $P\in S_0(\Omega)$.
\el

\bpf
Let $P\in S_0(\Omega)$ be a general point. 
Suppose $M_P$ is not a point. Then $M_P$ is a positive dimensional characteristic subspace. 
Choose a rank one boundary component $\Omega_\sigma$ such that $P\in \partial \Omega_\sigma$. Since $\Omega$ is non-tube type, 
$\Omega_\sigma$ is a ball of dimension at least $2$. Let $x\in \Omega_\sigma$. Choose a totally geodesic maximal tube type subdomain 
$\Omega_x\subset \Omega$ such that $x\in S_1(\Omega_x)$ and $P\in S_0(\Omega_x)$. 
Since $M_P$ is nontrivial,
by Corollary~\ref{image}, there exist nontrivial tube type subdomain $\Omega_1'(x)$ and a characteristic subdomain $\Omega_2(x)'$ such that 
$$f(\Omega_x)\subset \Omega_1'(x)\times \Omega_2'(x).$$
Moreover, in view of \eqref{def-LP}, $\Omega_2'(x)$ is completely determined by $Z_0(f(\mathscr S_P^1(\Omega_x))))$. 
Since $Z_0(f)$ is constant on each boundary component, $\Omega_2'(x)$ is parallel for all $x\in \Omega_\sigma$.
Since $\Omega_1'(x)\times \Omega_2'(x)$ is totally geodesic, $\Omega_1'(x)$ is also parallel for all $x\in \Omega_\sigma.$
Hence there exists a totally geodesic subspace of the form $\Omega_1'(\sigma)\times\Omega_2'(\sigma)\subset\Omega'$ of maximal rank such that $\Omega_1'(\sigma)$ is of tube type, $\Omega_2'(\sigma)$ is a characteristic subspace and 
$$f(\Omega_\sigma)\subset \partial(\Omega_1'(\sigma)\times\Omega_2'(\sigma)).$$
Since $\Omega_\sigma$ is a ball with dimension at least $2$, $f$ maps Shilov boundary to Shilov boundary and $\Omega_1'$ is of tube type, this implies that
$$f(\partial \Omega_\sigma)\subset\{A_P\}\times \partial\Omega_2'(\sigma)$$
for some $A_P\in S_0(\Omega_1')$ depending only on $f(P)$.
Since $\Omega_\sigma$ is arbitrary in $\mathscr S_P^1$ and $\Omega_2'(y)$ is parallel with $\Omega_2'(x)$ for $x, y\in \mathscr S_P^1(\Omega_x)$, $\Omega_2'(\sigma)$ is parallel with $\Omega_2'(\widetilde \sigma)$ for $\Omega_\sigma, \Omega_{\widetilde\sigma}\subset \mathscr S_P^1$ and 
$$f(\mathscr S_P)\subset\{A_P\}\times \Omega_2'.$$
Since $P$ is general, we obtain
$$f(S_1(\Omega))=f_1(\Omega)\times f_2(\Omega)\subset \Omega_1'\times\Omega_2'.$$
In particular, $f_1$ is constant. 
Since $f(S_0(\Omega))\subset S_0(\Omega')$, we obtain
$f_1(\Omega)$ is a point in $S_0(\Omega_1')$, contradicting the assumption that $f$ is proper.
\epf
\medskip

{\it Proof of Corollary~\ref{cor}}:
Let $f(S_0(\Omega)\cap U)\subset S_m(\Omega')$ for some open set $U$. If $m=0$, then by Theorem \ref{main}, $f$ is of the form
$$f=f_1\times f_2:\Omega\to \Omega_1'\times\Omega_2'\subset \Omega'$$
whose factor $f_1:\Omega\to\Omega_1'$ is a rational proper map.
If $\Omega_2'$ is trivial, then $f$ is rational. If $\Omega_2'$ is nontrivial, then ${\rm rank}(\Omega_1')<2q-1.$
Hence by Corollary 1.2 of \cite{K21}, $f_1$ has a standard embedding factor.
Suppose $m\geq 1.$ Then by the properness of $f$, after shringking $U$ if necessary,
$f(S_1(\Omega)\cap U)$ is contained in $S_{r}(\Omega')$ for some $r>m\geq 1$,
implying that
$$q'-r<2(q-1)$$
and on an open set of $S_1(\Omega)$, $f$ is transversal, i.e.,
$$f_*(\nu)\notin T^{1,0}_{f(x)}S_r(\Omega')+T^{0,1}_{f(x)}S_r(\Omega')$$
for all real vector $\nu\in T_xS_1(\Omega)$ transversal to $T_x^{1,0}S_1(\Omega)+T^{0,1}_xS_1(\Omega)$.
Then by Corollary 1.2 of \cite{K21}, $f$ has a standard embedding factor, which completes the proof.

\end{document}